\pdfoutput=1
\documentclass{aptpub}

\usepackage{hyperref}
\usepackage{fancyvrb}
\usepackage{graphicx}				% Use pdf, png, jpg, or eps§ with pdflatex; use eps in DVI mode
\usepackage{amssymb}

\usepackage{amsmath}
\usepackage{url}
\usepackage{enumerate}
\usepackage{mathtools}
\usepackage{algorithm}
\usepackage[table]{xcolor}% http://ctan.org/pkg/xcolor
\usepackage{makecell}
\usepackage{array}% http://ctan.org/pkg/array

\newcommand{\E}{\mathbb E}
\renewcommand{\P}{\mathbb P}

\newcommand{\argmin}{\mathop{\mathrm{argmin}}}

\newcommand{\Var}{\mathrm{Var}}

\authornames{GERGELY \'ODOR AND PATRICK THIRAN} % insert the authors here for use in running head. If three or more authors please use (for example) M.~YARROW {\it et al}. Author names should follow the same M.~YARROW format and if two authors, separate by 'AND'.
\shorttitle{Sequential metric dimension for random graphs} % insert short title here for use in running head

\begin{document}
\title{Sequential metric dimension for random graphs}

\authorone[EPFL]{Gergely \'Odor} 
\addressone{IINFCOM INDY2, \'Ecole Polytechnique F\'ed'erale de Lausanne, EPFL-IC Station 14 - B\^atiment BC, CH-1015 Lausanne, Switzerland} % Your postal address goes here.
\emailone{gergely.odor@epfl.ch} %Authors email goes here.

\authortwo[EPFL]{Patrick Thiran} 
\addresstwo{IINFCOM INDY2, \'Ecole Polytechnique F\'ed'erale de Lausanne, EPFL-IC Station 14 - B\^atiment BC, CH-1015 Lausanne, Switzerland} % Your postal address goes here.
\emailtwo{patrick.thiran@epfl.ch} %Authors email goes here.

\begin{abstract}

In the localization game on a graph, the goal is to find a fixed but unknown target node $v^\star$ with the least number of distance queries possible. In the $j^{th}$ step of the game, the player queries a single node $v_j$ and receives, as an answer to their query, the distance between the nodes $v_j$ and $v^\star$. The sequential metric dimension (SMD) is the minimal number of queries that the player needs to guess the target with absolute certainty, no matter where the target is.

The term SMD originates from the related notion of metric dimension (MD), which can be defined the same way as the SMD, except that the player's queries are non-adaptive. In this work, we extend the results of \cite{bollobas2012metric} on the MD of Erd\H{o}s-R\'enyi graphs to the SMD. We find that, in connected Erd\H{o}s-R\'enyi graphs, the MD and the SMD are a constant factor apart. For the lower bound we present a clean analysis by combining tools developed for the MD and a novel coupling argument. For the upper bound we show that a strategy that greedily minimizes the number of candidate targets in each step uses asymptotically optimal queries in Erd\H{o}s-R\'enyi graphs. Connections with source localization, binary search on graphs and the birthday problem are discussed.
\end{abstract}

%insert keywords separated by a semicolon. You should avoid including keywords which also appear in the title.

 % insert the primary Maths Subject Classification number in the first bracket
         % and the secondary ams number(s) in the second bracket
         % e.g. \ams{60E20}{49G03;49F10}
         %05C80 Random graphs (graph-theoretic aspects) 
         %68R05 Combinatorics in computer science
         %05C85 Graph algorithms (graph-theoretic aspects) 
         %60C05 Combinatorial probability

\section{Introduction}

With the appearance of new applications in network science, the theoretical analysis of search problems in random graph models is increasingly important. One such application is the source localization problem, where we assume that a stochastic diffusion process had spread over a graph starting from a single node, and we seek to find the identity of this node from limited observations of the diffusion process \cite{shah2011rumors, pinto2012locating}. If the diffusion models an epidemic outbreak, then our goal is to find patient zero, which is an important piece of information for both understanding and controlling the epidemic. The limited information about the diffusion is often the infection time of a small subset of sensor nodes \cite{pinto2012locating}. Recently, \cite{zejnilovic2013network} connected a deterministic version of the source localization problem with the metric dimension, a well-known notion in combinatorics introduced in 1975 by Slater \cite{slater1975leaves} and a year later by Harary and Melter~\cite{harary1976metric}.

\begin{defn}[$\mathrm{MD}$]
\label{MD_def}
Let $G=(V,E)$ be a simple connected graph, and let us denote by $d(v,w) \in \mathbb{N}$ the length of the shortest path between nodes $v$ and $w$. For $R=\{ w_1, \dots, w_{|R|}\} \subseteq V$ let $d(R,v) \in \mathbb{N}^{|R|}$ be the vector whose entries are defined by $d(R,v)_i=d(w_i,v)$. A subset $R\subseteq V$ is a \textit{resolving set}  in $G$ if $d(R,v_1)= d(R,v_2)$ holds only when $v_1= v_2$. The minimal cardinality of a resolving set is the \textit{metric dimension} (MD) of $G$.
\end{defn}

A resolving sensor set $R$ enables us to detect any epidemic source $v^\star$ when we assume that the observations are the vectors $d(R,v^\star)$, because every $v^\star$ generates a different deterministic observation vector. In many applications it is unrealistic to assume that $d(R,v^\star)$ can be observed, partially because we are not accounting for the noise in the disease propagation, but also because we infer $d(R,v^\star)$ from the time difference between the time the source gets infected and the time the sensors get infected, and usually we do not have access to the former information. We can define an analogous problem, where we do not assume that the time the infection began is known, if we require that all vectors in the set $\{ d(R,v^\star)+C \mid v^\star \in V, C\in \mathbb{Z} \}$ are different. The number of sensors needed in this scenario is called the double metric dimension (DMD) \cite{caceres2007metric}. Although the MD and the DMD can be very different in certain deterministic graph families, they seem to behave similarly in a large class of graphs, including Erd\H{o}s-R\'enyi random graphs \cite{spinelli2018many}, which are the focus of this paper. In this work, we consider only models that assume that the time the infection began is known. However, we believe our results can be extended to the DMD as well.

Previous work suggests that even with the assumption of deterministic distance observations, the number of sensors required to detect the source can be extremely large on real-word networks \cite{spinelli2017general}. One idea for mitigating this issue is to enable the sensors to be placed adaptively; using the information given by all previous sensors to place the subsequent ones \cite{zejnilovic2015sequential}. A substantial decrease in the number of required sensors in an adaptive version of the DMD compared to the DMD was observed experimentally by~\cite{spinelli2017general}. Intuitively, reducing the number of candidate nodes that could still be the source and focusing on only these candidate nodes can be very helpful, especially in real-word networks. However, it is not yet clear what property of the graph determines whether reduction in the number of required sensors is small or large. It is well known that in the barkochba game (binary search on a finite set), it does not matter whether the questions (queries) are non-adaptive or can be based on previous answers; the number of questions needed is $\lceil \log_2(N)\rceil$ in both cases. Source localization with non-adaptive (respectively, adaptive) sensor placement can be seen as a barkochba game with non-adaptive (respectively, adaptive) questions, where the questions are limited (to the nodes) and the answer does not have to be binary. It is the limitation on the available ``questions'' that creates a large gap between number of required ``questions'' in the adaptive and the non-adaptive version of the source localization problem. Our goal in this paper is to rigorously quantify this gap in source localization on a random graph model.

For the rigorous analysis, we consider an adaptive version of the MD in connected Erd\H{o}s-R\'enyi random graphs. In the combinatorics literature, this adaptive version of the MD was introduced in \cite{seager2013sequential} under the name sequential location number. The same notion was later referred to as the sequential metric dimension (SMD) by \cite{bensmail2018sequential}. We focus on the Erd\H{o}s-R\'enyi random graph model, because of the previous results on the MD of Erd\H{o}s-R\'enyi graphs by~\cite{bollobas2012metric}. The only other result on the MD of random graphs that we are aware of is the MD of uniform random trees \cite{mitsche2015limiting}. We do not consider this model in this paper, however, it is safe expect that the SMD would be significantly lower in this model than the MD. 

Some of the techniques used in this paper build directly on the techniques of \cite{bollobas2012metric}. The most important example is the \textit{expansion properties} of connected Erd\H{o}s-R\'enyi graphs, the main technique developed in \cite{bollobas2012metric}. According to this property, the observations $d(w_i,v)$ are dominated by one or two values from the set $\{D,D-1\}$, where $D$ is the diameter of $G$ (see Figure \ref{level_sets} and Table \ref{tools}). Hence, the information acquired in each step is essentially binary. In \cite{kim2015identifying}, the authors assume a very similar model to ours, except that the queries of the form $(v,r)$, and the answers are binary depending on whether the target is in the ball around node $v$ with radius $r$. Clearly, in Erd\H{o}s-R\'enyi graphs, where distance queries happen to have essentially binary answers, the two models are very similar. Indeed, \cite{kim2015identifying} independently recovers many of the results of \cite{bollobas2012metric}. In \cite{kim2015identifying}, the adaptive version of the problem is also intruduced, which is very similar to the SMD, but they do not have any results on the adaptive version of the problem in Erd\H{o}s-R\'enyi graphs. Since in the SMD we assume that we can use strictly more information than the binary model, our lower bounds are readily applicable to the binary model. The upper bounds are not readily applicable, but they could be extended with minimal modifications to the proof.

The binary nature of the answers to distance queries in Erd\H{o}s-R\'enyi graphs suggests that our setup has close connections with Generalized Binary Search \cite{nowak2011geometry}. In a sense, our problem setup can be seen as the dual version of graph binary search introduced by \cite{emamjomeh2016deterministic}, where the observations reveal the first edge in the shortest path instead of its length. Although the two models share some similarities, we must point out that while \cite{emamjomeh2016deterministic} focuses on an algorithm for general graphs (with noisy but adaptive observations), our work provides asymptotically almost sure results on the sample complexity of an algorithm and a matching lower bound for all possible algorithms on Erd\H{o}s-R\'enyi graphs (with noiseless observations), thus the authors aim for different goals. In terms of goals, the work most similar to ours is perhaps \cite{dudek2019note}. In their paper, the authors consider a version of the Cop and Robber game on Erd\H{o}s-R\'enyi graphs of diameter two: The target can ``move'' between turns, and in order to locate this moving target, it is not the number of turns but the number of sensors that the player selects in each turn that we want to minimize. Recently, the results of \cite{dudek2019note} were extended by \cite{dudek2019localization} to Erd\H{o}s-R\'enyi graphs with a diameter larger than two, and they found that in that range, the number of sensors needed in the Cop and Robber game is strictly less then the SMD. Similarly to our proofs, the proofs in \cite{dudek2019localization} make use of the expansion properties developed in \cite{bollobas2012metric}. Finally, we also mention the recent work by \cite{lecomte2020noisy}, where a noisy version of the SMD is studied in path graphs.

The methods in this paper connect several different ideas developed in different communities; these ideas have not been connected before. In Section \ref{warmup1}, we abstract out one of the key ideas of \cite{bollobas2012metric} and connect it with the Birthday Problem. In Section \ref{warmup}, we connect the SMD with Generalized Binary Search \cite{nowak2011geometry}. In Section \ref{metric}, we introduce the expansion properties of $\mathcal{G}(N,p)$ random graphs. In Section \ref{main_result}, we combine the ideas of Sections \ref{warmup} and \ref{metric} and present the main result of this paper. Finally, we conclude the paper with a discussion in Section \ref{discussion}, and we place some of the proofs to the end of the paper for better readability.

         % Introductory Chapter
\section{Problem Statement and Summary of Results}
\label{definitions}

\subsection{Problem Statement}

Although we explained our problem with the source localization problem, in the rest of the paper we adopt the vocabulary of binary search. Let $v^\star \in V$ be the \textit{target} node. The target node is unknown to us, but for a set \textit{queries} $R \subseteq V$ the distance $d(R,v^\star)$ is known.

\begin{defn}[candidate targets]
\label{def:candi}
Given a set of queries $R$, the set of \textit{candidate targets} for the graph $G$ is 
$$\mathcal{T}_{R}(G) = \{ v \in V \mid d(R,v)=d(R,v^\star) \}.$$
\end{defn}

Our goal is to detect $v^\star$, which means we would like a set $R$ with $\mathcal{T}_{R}(G)=\{v^\star\}$, or equivalently $|\mathcal{T}_{R}(G)|=1$ (as $v^\star \in \mathcal{T}_{R}(G)$ must always hold). Recall, that for a resolving set $R$ we have $|\mathcal{T}_{R}(G)|=1$ for every $v^\star \in V$. In contrast, in the adaptive case, a (potentially) different $R$ is constructed for every $v^\star$; in the $j^{th}$ step we select query $w_j$ based on the distance information revealed by $R_{j-1}=\bigcup_{k=1}^{j-1} w_k$, and we still aim for $|\mathcal{T}_{R_{j}}(G)|=1$.

\begin{defn}[$\mathrm{SMD}$]
\label{def_smd}
Let $\mathrm{ALG}(G)$ be the set of functions $$g: \{(G,R,d(R,v^\star)) \mid R \subseteq V, v^\star\in V\} \rightarrow V.$$ The \textit{sequential metric dimension} (SMD) of $G$ is the minimum $r \in \mathbb{N}$ such that there is a query selection algorithm $g \in \mathrm{ALG}(G)$, for which if we let $R_0=\emptyset$ and $R_{j+1} = R_j \cup g(G,R_{j}, d(R_{j},v^\star))$, then $|\mathcal{T}_{R_{r}}(G)|=1$ for any $v^\star \in V$.
\end{defn}

The set $\mathrm{ALG}(G)$ might contain functions that are not computable in polynomial time, hence we define a slightly stricter notion of SMD where the next query has to be polynomial time computable (SMDP).

\begin{defn}[$\mathrm{SMDP}$]
Let $\mathrm{PALG}(G)$ be the subset of $\mathrm{ALG}(G)$ with polynomial time complexity. Then, SMDP(G) is the minimum $r \in \mathbb{N}$ such that there is a query selection algorithm $g \in \mathrm{PALG}(G)$, for which if we let $R_0=\emptyset$ and $R_{j+1} = R_j \cup g(G,R_{j}, d(R_{j},v^\star))$, then $|\mathcal{T}_{R_{r}}(G)|=1$ for any $v^\star \in V$.
\end{defn}

\begin{figure}[h]
\begin{center}
  \includegraphics[width=0.7\textwidth]{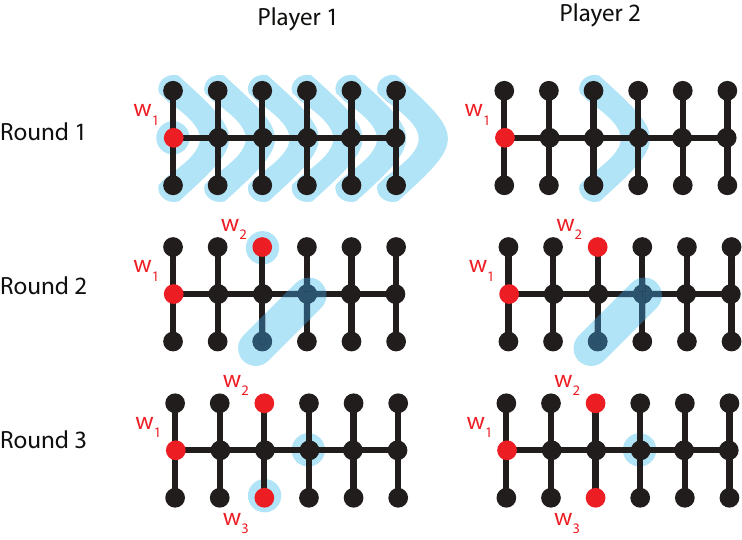}
  \caption{An example of how the SMD can be interpreted as a two-player game. In the $j^{th}$ round, Player 1 creates the set $R_{j}$ by adding a sensor node $w_j$ (marked in red) to $R_{j-1}$. The sensor $w_j$ partitions the current candidate target set $\mathcal{T}_{R_{j-1}}(G)$ based on distances (marked in blue). In turn, Player 2 must provide a distance from $w_j$ to a feasible but not necessarily predetermined source node, which is equivalent to selecting one of the blue sets. Player 1 tries to reduce and Player 2 tires to increase the total number of rounds until the end of the game, which happens when $\mathcal{T}_{R_{j}}(G)$ shrinks to a single element. In this example, the game ends in 3 rounds if both players play optimally. Hence, the SMD of this ``comb graph'' of size $18$ is also $3$. In fact, the SMD of the ``comb graph'' of any size $n\ge9$ is still 3, in sharp contrast with the MD of the same graph, which is $n/3$. } 
  \label{SMD_example}
  \end{center}
\end{figure}

The definition of SMD and SMDP is intrinsically algorithmic. It is useful to think of them as two-player games. In each step, Player 1 selects a query and tries to reduce as fast as possible the candidate set to a single element. Player~2 must then provide an observation that is consistent with at least one of the target nodes. If there are multiple such observations, Player 2 can choose one to try to make the game as long as possible. In this setting Player 2 does not decide on the source $v^\star$ in advance, but must always be consistent with the observations that have been revealed so far (i.e., $\mathcal{T}_{R_j}(G)$ can never be empty). Since every predetermined source $v^\star$ can be found this way by Player 1, and since for every set of answers provided by Player 2 there is a node that could have been the source, the SMD can be seen as the number of steps the game takes if both players play optimally, and SMDP is the same if Player 1 must compute their next move in polynomial time in each step. See Figure \ref{SMD_example} for an example of how the two-player game corresponding to the SMD is played.

Clearly $1\le \mathrm{SMD}\le \mathrm{SMDP} \le \mathrm{MD} \le N$, as being able to adaptively select the queries only gives Player 1 more power. Before we proceed to compute the difference between in Erd\H{o}s-R\'enyi graphs, we first introduce two easier problems defined on matrices as warmups and we address the problem on graphs in Section \ref{main_result}.

\subsection{Summary of Results}

To ease notation in this short summary, we express our main results on the SMD in terms of the MD. The precise formulation and the proof of our main theorem is in Section \ref{main_result}. In its crudest form, our main result says that the ratio of the SMD and the MD is between $1$ and $1/2$ a.a.s. From this statement and the results of \cite{bollobas2012metric}, it is already possible to infer the asymptotic behavior of the SMD. In \cite{bollobas2012metric} it was found that $\mathrm{MD}(\mathcal{G}(N,N^{x-1}))=N^{1-\lfloor 1/x \rfloor x+o(1)}$, which means that the MD is a (non-monotonically changing, ``zig-zag'' shaped) power of $N$, unless $1/x$ is an integer, in which case the MD is a constant times $\log(N)$. In this paper we prove that the same is true for the SMD, hence the crude asymptotic behaviour of the SMD is completely characterized. 

Our results enable us to make a more precise statement on the ratio of the SMD and the MD. For the values of $p$ where the MD is logarithmic, we are able to exactly determine the leading constant of the ratio. Since the dependence of this leading constant on $p$ and $N$ is rather complicated, we simply denote it by $F_\gamma(p,N)$ in Theorem~\ref{meta_thm}, and we defer the precise definition of $F_\gamma(p,N)$ to Remark~\ref{rem:Fgammaeta}. It is shown in Remark~\ref{rem:Fgammaeta} that $F_\gamma(p,N)$ can take any value in the interval $(1/2,1)$, which implies that there are values of $p$ for which the SMD is strictly smaller than the MD.

For the values of $p$ where the MD is a power of $N$, we have a lower bound on the ratio of the SMD and the MD, which in general does not match the trivial upper bound given by the observation $\mathrm{SMD}\le \mathrm{MD}$. In Theorem \ref{meta_thm}, we denote our lower bound by $F_\eta(p,N)$, the definition of which we again defer to Remark~\ref{rem:Fgammaeta}, but we mention that $F_\eta(p,N)$ takes values in the interval $(1/2,1)$, which means that there is a non-zero $1-F_\eta(p,N)$ gap between our upper and lower bounds for these values of $p$. We conjecture that this gap is due to the use of the first moment method in our proofs, and that the lower bound can be improved to $1$ with more sophisticated techniques.

\begin{thm}
\label{meta_thm}
Let $N \in \mathbb{N}$ and $p \in [0,1]$ such that $\frac{\log^5(N)}{N}\ll p$ and $\frac{1}{\sqrt{N}} \ll 1-p$. Let $G$ be a realization of a $\mathcal{G}(N,p)$ random graph. Then,
\begin{align}
1\ge \frac{\mathrm{SMD}(G)}{\mathrm{MD}(G)} &= F_\gamma(p,N) +o(1) \ge \frac12 + o(1) \quad \text{if } (Np)^i=\Theta(N) \text{ for } i \in \mathbb{N}\\
1\ge \frac{\mathrm{SMD}(G)}{\mathrm{MD}(G)} &\ge F_\eta(p,N) +o(1) \ge \frac12  +o (1) \quad \text{otherwise},
\end{align}
hold a.a.s, where $F_\gamma$ and $F_\eta$ are functions of $p, N$ that are explicitly expressed in Remark~\ref{rem:Fgammaeta}.
\end{thm}

\begin{figure}[h]
\begin{center}
  \includegraphics[width=\textwidth]{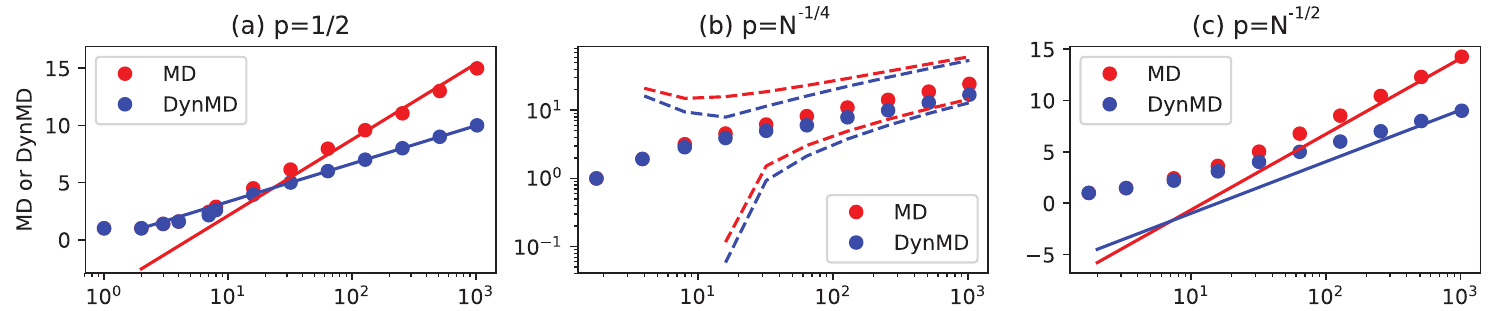}
  \caption{The red and blue dots show the approximated value of the MD and the SMD of simulated Erd\H{o}s-R\'enyi graphs computed by the toolbox \cite{odorgergo_2020_3937453} averaged over 100 iterations (confidence intervals are too small to be plotted). The slope of the red and blue lines is computed by Theorem \ref{main_thm} and the intercept is chosen to fit the last few data points. On (semi-log) plots (a) and (c) we have $(Np)^i=\Theta(N)$ for $i=0$ and $i=1$ respectively. For such parameters the MD and the SMD are both logarithmic and there is a constant factor difference between them. On the contrary, on (log-log) plot (b) we have $(Np)^i \ne \Theta(N)$ for all $i \in \mathbb{N}$. For such parameters the MD and the SMD grow as a power of $N$, and there is a gap between the theoretical upper and lower bounds, which as shown with dashed curves.}
  \label{MD-SMD-assymptotics}
  \end{center}
\end{figure}

See Figure \ref{MD-SMD-assymptotics} for simulation results confirming Theorem \ref{meta_thm}, and its more precise version, Theorem \ref{main_thm}.

In Theorem \ref{meta_thm}, we were able to succinctly state our main results in terms of the MD, however, in our proofs we cannot take such shortcuts. Instead of directly using the results of \cite{bollobas2012metric} on the MD, we use some of their techniques, and we complement them with new techniques of our own. For instance in the SMD upper bound, we need to analyse an interactive game of possibly $N$ steps instead of selecting the queries in a single round. In particular, the order in which we reveal the edges of the random graph is completely different in our analysis than in \cite{bollobas2012metric}. For the SMD lower bound, we could have split our proof in several cases, and for some of them we could have used the results of \cite{bollobas2012metric} directly. Instead, we introduce a coupling argument, which succeeds without case-work, and gives a clean alternative proof to the MD lower bound as well. 

\section{Warmup1: Random Bernoulli Matrices with Pairwise Different Columns}
\label{warmup1}

In this section we consider an $M \times N$ random matrix $A$, with entries drawn independently from a Bernoulli distribution, and we are interested in the minimal $M$ for which $A$ still has pairwise different columns with high probablity. This $M$ can be viewed as the query complexity of binary search with random Bernoulli queries, where the $i^{th}$ query can distinguish between targets $j$ and $k$ if $A_{ij}\ne A_{ik}$. Anothers

For notation, let us consider the binary matrix $A$ with row indices $\mathcal{R} = [M]$ and column indices $\mathcal{C} = [N]$. For $R\subseteq \mathcal{R}$ and $W \subseteq \mathcal{C}$, let $A_{R,W}$ be the submatrix of $A$ restricted to rows $R$ and columns $W$.

\begin{thm}
\label{warmup1_thm}
Let $N\in \mathbb{N}$, let $0 < q(N) \le 1/2$ and $M(N) \in \mathbb{N}$ be functions possibly depending on $N$, and let us define the random matrix $A \in \mathrm{Ber}(q)^{M\times N}$. Let $\mathcal{A}$ be the property that $A$ has pairwise different columns. Then 
\begin{equation}
\hat{M}(N)= \frac{\log(N)}{\log\left(1/\sqrt{q^2+(1-q)^2}\right)}
\end{equation}
 is the threshold function for $\mathcal{A}$ . That is, for any $0 < q(N) \le 1/2$ and $1 \gg \epsilon(N)\gg~\frac{1}{\log(N)}$,

\begin{enumerate}[(i)]
\item if $M\ge (1+\epsilon(N))\hat{M}$, then $\lim\limits_{N\rightarrow \infty}P(A\in \mathcal{A})=1$
\item  if $M\le (1-\epsilon(N))\hat{M}$, then $\lim\limits_{N\rightarrow \infty}P(A\in \mathcal{A})=0$.
\end{enumerate}
\end{thm}

We could not find this particular theorem stated this way in the literature, however, there exist many related results. Computing the probability that an $N \times N$ random Bernoulli matrix is singular is a famous problem first proposed by Koml\'os in 1967 \cite{komlos1967determinant}. Clearly, if the matrix has two identical columns, then it is also singular, hence we obtain the lower bound

$$\P(A\not \in \mathcal{A}) \le \P(A \text{ is singular}).$$

Most of the research on the singularity of random Bernoulli matrices has been on the upper bound \cite{kahn1995probability,tao2006random}, with the exception of \cite{arratia2013singularity}. In \cite{arratia2013singularity}, the authors lower bound $\P(A\not \in \mathcal{A})$ by using an inclusion-exclusion type argument. However, this bound is too loose in our case as we are interested in $P(A\in \mathcal{A})$ of an $M \times N$ matrix, where $M$ is close to the threshold. Our analysis in this paper could potentially be applied to tighten some of the bounds in \cite{arratia2013singularity}, although the improvement would appear only in a high ($5^{th}$) order term of the bound.

Another well-studied problem related to $P(A\in \mathcal{A})$ is the Birthday Problem (BP). Indeed, when $q=1/2$, we obtain the standard formulation of the BP with $N$ people and $2^M$ days. For $0<q<1/2$, the columns (birthdays) are not equiprobable anymore, hence we obtain BP with heterogenous birthday probabilities. The non-coincidence probability of two birthdays in this case has been computed exactly using a recursive formula by \cite{mase1992approximations}. Rigorous closed-form approximations for constant $q$s were given by \cite{arratia1989two}. Intuitively, the events that two birthdays coincide are rare and almost independent, so the number of coincidences can be approximated by a Poisson random variable. However, Poisson approximation can work only as long as the number of pairwise collisions dominates the number of multi-collisions (i.e., collisions of $\ge 3$ columns), which happens only for $q$s that do not decrease too fast with $N$. For fast-decaying $q$s, we need to use a different technique; we upper bound $P(A\in \mathcal{A})$ by the probability that $A$ does not contain two identically zero columns. Indeed, the event that all columns are different implies that that no two identically zero columns can exist, hence it must have a smaller probability.

Theorem \ref{warmup1_thm} can also be viewed as a simplified version of \cite{bollobas2012metric}, which will become clear in Section \ref{main_result}. Our proofs also follow \cite{bollobas2012metric}: We chose to use a combination of the first moment method and Suen's inequality \cite{suen1990correlation,Janson1998NewVO} is used instead of Poisson approximation. 

\subsection{Proof of Theorem \ref{warmup1_thm} Part (i) }
\label{p1ofT1}
\begin{proof}
Let $M=\frac{(1+\epsilon(N))\log(N)}{\log\left(1/\sqrt{q^2+(1-q)^2}\right)}$ with $\epsilon(N)\gg \frac{1}{\log(N)}$ and $X$ be the number of pairs of columns in $\mathcal{C}$ which are identical (i.e., colliding). Let $X_{xy}$ be the indicator of  $A_{\mathcal{R},x}=A_{\mathcal{R},y}$ for $x\ne y \in \mathcal{C}$, and let $P_k$ be the marginal that $k$ fixed columns all collide. By the identity 
\begin{equation}
\label{identity}
\alpha^{\hat{M}}=\alpha^{\frac{\log(N)}{\log\left(1/\sqrt{q^2+(1-q)^2}\right)}}=N^{-\frac{2\log(\alpha)}{\log\left(q^2+(1-q)^2\right)}}
\end{equation}
for all $\alpha \in \mathbb{R}$, we have by taking $\alpha=(q^2+(1-q)^2)$

\begin{align}
\E[X]&=\E[\sum\limits_{x \ne y \in \mathcal{C}} X_{xy}] =\sum\limits_{x\ne y \in \mathcal{C}}P_2=\binom{N}{2} (q^2+(1-q)^2)^{(1+\epsilon(N))\hat{M}}  \nonumber \\
&=\frac{N(N-1)}{2}N^{-2(1+\epsilon(N))} < \frac12N^{-2\epsilon(N)} \rightarrow 0 
\end{align}

as $N\rightarrow\infty$. By the first moment method, it follows as $N\rightarrow\infty$ that

\begin{equation}
\P(A \not\in \mathcal{A}) = \P(X>0) \le \E[X] \rightarrow 0
\end{equation}
\end{proof}

\subsection{Proof of Theorem \ref{warmup1_thm} Part (ii)}
\label{p2ofT1}
The lower bound will be more involved. We are going to split our argument into two cases; Case 1: $q> \epsilon$ and Case 2: $q \le \epsilon$. In Case 1, we will use Suen's inequality and in Case 2 we will bound the probability that $A$ does not contain two identically zero columns.

\begin{proof}[Proof of Theorem \ref{warmup1_thm} Part (ii), Case 1: $q> \epsilon$]

%First, let us fix $\epsilon=\Theta(\frac{1}{\log(N)})$. It is enough if we show the statement for this particular $\epsilon$, since property $\mathcal{A}$ is monotone in $M$. Clearly, adding more columns can only decrease the collision probability (we can see this by e.g. two round exposure). 
To be able to apply Suen's inequality, we will need to show that pairwise collisions dominate three-way collisions. For this we must estimate $P_2$ and $P_3$. Using equation \eqref{identity}, with $\alpha=(q^2+(1-q)^2)^{(1-\epsilon)}$

\begin{equation}
\label{P2}
P_2=(q^2+(1-q)^2)^{(1-\epsilon)\hat{M}} =N^{-2(1-\epsilon)},
\end{equation}

and with $\alpha=(q^3+(1-q)^3)^{(1-\epsilon)}$

\begin{equation}
\label{P3}
P_3=(q^3+(1-q)^3)^{(1-\epsilon)\hat{M}}=N^{-2(1-\epsilon)\frac{\log(q^3+(1-q)^3)}{\log\left(q^2+(1-q)^2\right)}}\le N^{-(1-\epsilon)(3+\frac32 q)},
\end{equation}

because one can check that for $0 \le q\le \frac12$

\begin{equation}
\frac{\log(q^3+(1-q)^3)}{\log\left({q^2+(1-q)^2}\right)}\ge \frac{3}{2}+\frac{3q}{4}.
\end{equation}

Now we apply Suen's inequality \cite{suen1990correlation,Janson1998NewVO} to our setting. Let us define the index set $I= \{ \{x,y\} \mid x\ne y \in \mathcal{C}\}$, which allows us to index $X_{xy}$ as $X_{\alpha}$ for $\alpha \in I$. For each $\alpha \in I$ we define the ``neighborhood of dependence'' $B_\alpha=\{ \beta \in I \mid \beta \cap \alpha \ne \emptyset \}$. Indeed, $X_\alpha$ independent of $X_\beta$ if $\beta\not\in B_\alpha$. Then, Suen's inequality implies

\begin{equation}
\label{suen}
\P(A \in \mathcal{A})=\P(X=0) \le \exp(-\lambda+\Delta e^{2\delta})
\end{equation}

where, using $|I|=\binom{N}{2}$, $|B_\alpha|=(2N-3)$ and equations \eqref{P2} and \eqref{P3},

\begin{align}
\lambda&=\sum\limits_{\alpha \in I} \E[X_\alpha] = \binom{N}{2} P_2 > \frac14 N^{2\epsilon}  \\
\Delta&=\sum\limits_{\alpha \in I} \sum\limits_{\alpha \ne \beta \in B_\alpha} \frac12 E[X_{\alpha}X_{\beta}] =  \binom{N}{2}(2N-4) \frac12 P_3 \le N^{3-(1-\epsilon)(3+\frac32 q)} \\
\delta&=\max\limits_{\alpha \in I} \sum\limits_{\alpha \ne \beta \in B_\alpha} E[X_{\beta}]= (2N-4)P_2 < 
2N^{-1+2\epsilon}.
\end{align}

We note that as $N \rightarrow \infty$ we have $1<e^{2\delta}<e^{4N^{-1+2\epsilon}}\rightarrow 1$. Hence, for $N$ large enough we have

\begin{equation}
\label{lamdel}
-\lambda + \Delta e^{2\delta} <  -\lambda + 2\Delta <  -\frac14 N^{2\epsilon} + 2N^{3-(1-\epsilon)(3+\frac{3}{2}q)}= N^{2\epsilon}\left(- \frac14+2N^{\epsilon - \frac{3}{2}q(1-\epsilon)} \right) \rightarrow -\infty
\end{equation}

because when $q > \epsilon$ and $\frac{3(1-\epsilon)}{2}>1$ we have $N^{\epsilon - \frac{3}{2}q(1-\epsilon)} \rightarrow 0$. By Equation \eqref{suen} we can conclude $\P(A \in \mathcal{A}) \rightarrow 0$.

We see that for smaller $q$s such a Suen's inequality type analysis cannot work because the number of colliding triples ($\Delta$), starts dominating the number of colliding pairs ($\lambda$).
\end{proof}

\begin{proof}[Proof of Theorem \ref{warmup1_thm} Part (ii), Case 2: {$q\le \epsilon$}]

We would like to upper bound the probability that all columns are distinct by the probability that at most one column is identically $0$. If we denote the number of identically $0$ columns by $Z$, we want to show that $P(Z<2)\rightarrow 0$. We start by proving $\E[Z] \rightarrow \infty$.

One can check that for $0 \le q \le \frac14$ (which we may assume as $q\le \epsilon \rightarrow 0$),
\begin{equation}
\label{fact}
\frac{\log(1-q)}{\log(\sqrt{q^2+(1-q)^2)}} \le 1+\frac{9}{10}q.
\end{equation}

Then, first using equation \eqref{identity} with $M=(1-\epsilon)\hat{M}$ and $\alpha=1-q$, and next applying inequality \eqref{fact} and the assumptions $q\le \epsilon$ and $\epsilon\gg \frac{1}{\log(N)}$, we have

\begin{align}
\E[Z]&=N(1-q)^M=N^{1+(1-\epsilon)\frac{\log(1-q)}{\log(1/\sqrt{q^2+(1-q)^2)}}} \ge N^{1-(1-\epsilon)(1+\frac{9}{10}q)} \nonumber \\
&= N^{ \epsilon -\frac{9}{10}q(1-\epsilon)} \ge N^{ \epsilon -\frac{9}{10}\epsilon} \rightarrow \infty
\end{align}

We are going to use a standard concentration bound to finish this proof.

\begin{lem}[Chernoff bound]
\label{chernoff}
Let X be a binomial random variable. Then, for $0<\tau <1$ we have
$$P(|X - \E[X]| \ge \tau E[X ]) \le 2 e^{-\frac{\tau^2 \E[X]}{3}}.$$
\end{lem}

By Lemma \ref{chernoff}, and since for $N$ large enough we have $\E[Z]>2$,

\begin{align}
\label{EZ}
\P(A \in \mathcal{A}) & \le \P(Z< 2)=\P\left(Z-\E[Z]<  2-\E[Z]\right) \nonumber \\
&\le \P\left(|Z-\E[Z]|\ge \left(1-\frac{2}{\E[Z]}\right)\E[Z]\right) \nonumber \\
 & \le 2e^{-\frac13\left(1-\frac{2}{\E[Z]}\right)^2\E[Z]} \le 2e^{-\frac13\left(1-\frac{4}{\E[Z]}\right)\E[Z]}= 2e^{\frac43} e^{-\frac{E[Z]}{3}} \rightarrow 0
\end{align}
\end{proof}

%%%%%%%%%%%%%%%%%%%%%%%%%%%%%%%%%%% WARMUP 2

\section{Warmup2: Identifying Codes or Binary Search with Randomly Restricted Queries}
\label{warmup}

In the previous section, we treated binary search with completely random entries. In this section, the queries will be selected by us, but we may only choose from a random subset of all queries (of size $N$). We start by adapting the definitions in Section \ref{definitions} to the matrix setup used in Section \ref{warmup1}.

\begin{defn}[QC]
Let $A \in \{0,1\}^{N\times N}$ be a binary matrix. Let the \textit{query complexity} ($\mathrm{QC}$) of $A$ be the minimum $|R|$ such that $A_{R,\mathcal{C}}$ has pairwise different columns.
%$$\min\{|R| \mid A_{R,[N]} \text{ has pairwise different columns}\}.$$
%the minimum $j^\star$ such that there is $R\subseteq [N]$ of cardinality $j^\star$, such that in $A_{R,[N]}$ any two columns have at least one different element.
\end{defn}

If the submatrix $A_{R,\mathcal{C}}$ has pairwise different columns, the set $R$ is also called an \textit{identifying code} \cite{karpovsky1998new} of the graph which has adjacency matrix $A$ (we must also allow self-edges to have equivalence, which is usually not part of the definition). Identifying codes are are closely related to resolving sets. If the graph has diameter two, the only difference between the two notions is that in the case of resolving sets we may receive three kinds of measurements (0, 1 and 2), not just two. However, we receive the 0 measurement only if we accidentally query the target, which can be ignored in many cases, hence the information we get is essentially binary for resolving sets as well. 

Identifying codes of Erd\H{o}s-R\'enyi graphs have been studied by \cite{frieze2007codes}. In fact, \cite{frieze2007codes} already featured some of the ideas that lead to characterizing the MD of Erd\H{o}s-R\'enyi graphs in~\cite{bollobas2012metric}. Part of the main theorem in this section (on the QC of Bernoulli random matrices, Theorem \ref{warmup2_thm}) has also appeared in \cite{frieze2007codes} for a limited range of parameters, which we extend using the tools of \cite{bollobas2012metric}. Our proof of the QC of Bernoulli random matrices does not feature any new ideas, we only include it for the sake of completeness. We also define the adaptive version of the problem, the \textit{sequential query complexity} (SQC), which will be similar to the $\mathrm{SMD}$. The upper bound on the SQC will be algorithmic and will be quite different from the tools in \cite{frieze2007codes} and \cite{bollobas2012metric}.

\begin{defn}[candidate targets]
\label{def:candi_A}
Given a set of queries $R$ and target $v^\star$, the set of \textit{candidate targets} for the matrix $A$ is 
$$\mathcal{T}_{R}(A) = \{ v \in [N] \mid A_{R,\{v\}} = A_{R,\{v^\star\}} \}.$$
\end{defn}

\begin{defn}[SQC and SQCP]
Let $\mathrm{ALG}(G)$ (respectively, $\mathrm{ALGP}(G)$) be the set of functions (respectively, polynomial time computable frunctions) $$g: \{(A,R,A_{R,\{v^\star\}}) \mid R \subseteq \mathcal{R}, v^\star\in \mathcal{C}\} \rightarrow \mathcal{R}.$$ The \textit{sequential query complexity} $\mathrm{SQC}(G)$ (respectively, $\mathrm{SQCP}(G)$) is the minimum $r \in \mathbb{N}$ such that there is a query selection algorithm $g \in \mathrm{ALG}(G)$ (respectively,  $g \in \mathrm{ALGP}(G)$), for which if we let $R_0=\emptyset$ and $R_{j+1} = R_j \cup g(G,R_{j}, A_{R_{j-1},\{v^\star\}}))$, then $|\mathcal{T}_{R_{r}}(A)|=1$ for any $v^\star \in \mathcal{C}$.
\end{defn}

\begin{figure}[h]
  \label{QC-SQC}

\begin{center}
  \includegraphics[width=0.7\textwidth]{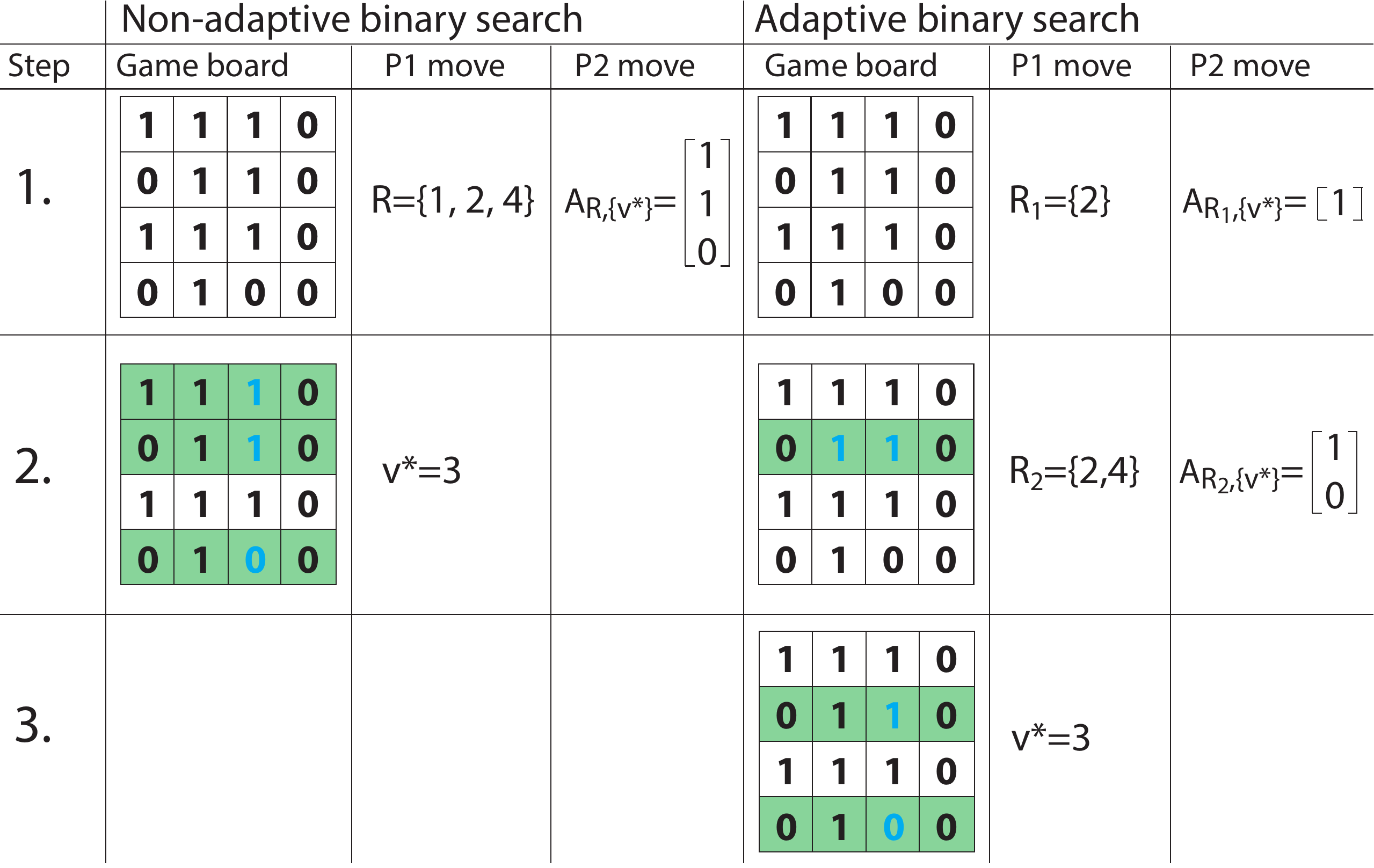}
  \caption{The process of non-adaptive and adaptive binary search with restricted queries with target $v^\star=3$. The queries are marked with green and the observations are marked with blue. }
  \label{QC}
\end{center}
\end{figure}

We show the difference between the  $\mathrm{QC}$ and $\mathrm{SQC}$ on an example in Figure \ref{QC}. The advantage studying the $\mathrm{QC}$ and $\mathrm{SQC}$ before the $\mathrm{MD}$ and $\mathrm{SMD}$ is that we have simpler results without the small dependencies that always arise in the graph setting.

%%%%%%%%%%%%%%%%%%%%%%%%%%%%%%%%%%%% THEOREM
\begin{thm}
\label{warmup2_thm}
Let $N \in \mathbb{N}$, let $0< q <1$, let $A \in \mathrm{Ber}(q)^{N\times N}$ and $\gamma_{sqc}=\max(q,1-q)$. Then a.a.s,

\begin{enumerate}[(i)]
%\item If $1-\gamma_{sqc}=\Theta(1)$ then we have
 %$$\mathrm{SQC}(A)= (1+o(1))\frac{\log(N)}{\log(1/\gamma_{sqc})}=\Theta(\log(N))$$
\item If $1\ge 1-\gamma_{sqc}\gg \frac{\log(N)}{N}$ then with $\eta=1+\log_N(\log(1/\gamma_{sqc}))$ we have
\begin{align}
& \mathrm{SQC}(A) \ge (\eta-o(1))\frac{\log(N)}{\log(1/\gamma_{sqc})} \nonumber \\ 
& \mathrm{SQCP}(A)  \le (1+o(1))\frac{\log(N)}{\log(1/\gamma_{sqc})}.
\end{align}
\item If $\frac{\log(N)}{N} \gg 1-\gamma_{sqc}$ then $\mathrm{SQC}(A)$ is undefined.
\end{enumerate}

The results for $\mathrm{QC}(A)$ are of the same form, except that instead of $\gamma_{sqc}$ we have $\gamma_{qc}=\sqrt{q^2+(1-q)^2}$. 
\end{thm}

\begin{rem}
\label{rem1}
If $1-\gamma_{sqc}=o(1)$ (which is equivalent to $1-\gamma_{qc}=o(1)$), we have
\begin{equation}
\label{rem1_eq}
\frac{1}{\log(1/\gamma_{qc})}=(1+o(1))\frac{1}{\log(1/\gamma_{sqc})}= (1+o(1))\frac{1}{1-\gamma_{sqc}},
\end{equation}
In this case $\mathrm{SQC}(A)=(1+o(1))\mathrm{QC}(A)=\omega(\log(N))$, so the SQC and the QC have the same asymptotic behaviour.
\end{rem}

\begin{rem}
\label{rem2}
If $1-\gamma_{sqc}=\Theta(1)$, then $\eta \rightarrow 1$, so the upper and lower bounds match in part (i) up to a multiplicative factor tending to one. In this case $\mathrm{SQC}(A)=\Theta(\log(N))$ and  $\mathrm{QC}(A)=\Theta(\log(N))$. 
\end{rem}

\subsection{Connection between Theorems \ref{warmup1_thm} and \ref{warmup2_thm} }
\label{interpret_1_and_2}
The main difference between the two theorems is that in Theorem \ref{warmup1_thm} we sample $M$ queries and use all of them, whereas in Theorem \ref{warmup2_thm} we sample $N$ queries and select (adaptively or non-adaptively) only a subset of them. The subset we select is of size $\mathrm{SQC}(A)$ or $\mathrm{QC}(A)$. Of course if $\gamma_{sqc}$ is so close to one that even all of the $N$ queries are not sufficient to locate the target, then it is impossible to find the target in the adaptive and the non-adaptive case as well. This intuition is made more precise in the following remark.

\begin{rem}
\label{thm1_implies_p3}
If $\frac{\log(N)}{N} \gg 1-\gamma_{sqc}$ the lower bound in part (i) would give $\mathrm{SQC}>N$, which is a contradiction, hence part (ii) can also be viewed as a special case part (i). Moreover, part (ii) is also implied by Theorem \ref{warmup1_thm}. Indeed, suppose $\frac{\log(N)}{N} \gg 1-\gamma_{sqc}$. After reordering

\begin{equation}
N \ll \frac{\log(N)}{1-\gamma_{sqc}} = (1+o(1))\frac{\log(N)}{\log(1/\gamma_{qc})} 
\end{equation}

because of equation \eqref{rem1_eq}. Then, by Theorem \ref{warmup1_thm} with $q'=1-\gamma_{sqc} \in (0, 1/2]$, we know that $A'=\mathrm{Ber}(q')^{N\times N}$ has two identical columns a.a.s. This implies that $A=\mathrm{Ber}(q)^{N\times N}$ has two identical columns, and thus $\mathrm{SQC}(A)$ and $\mathrm{QC}(A)$ are undefined a.a.s. Similarly, if $\frac{\log(N)}{N} \ll 1-\gamma_{sqc}$, then $A$ does not have two identical columns, which implies that $\mathrm{SQC}(A)$ and $\mathrm{QC}(A)$ are well-defined.
\end{rem}

The proofs of Theorems \ref{warmup1_thm} and \ref{warmup2_thm} are quite similar, except for the $\mathrm{SQC}$ upper bound. However, while in Theorem \ref{warmup1_thm} we have matching upper and lower bounds, in Theorem \ref{warmup2_thm} there is an $1-\eta$ gap between them. It is an open question whether this gap can be closed, and we conjecture that the constant in front of the lower bound can be improved from $\eta$ to $1$.

Let us also give intuition about the new notation. On the range $q\in (0,\frac{1}{2}]$, the variable $1-\gamma_{sqc}$ is just $q$, and it serves essentially the same purpose. The reason for introducing a new variable is that we can highlight the symmetry of the adaptive and non-adaptive case here and later in the text. The other new variable is $\eta$, which is monotonically decreasing in $\gamma_{sqc}$. We note that since $\gamma_{sqc}\ge \frac12$, we always have $\eta<1$, hence the upper and lower bounds never contradict. However, $\eta$ can be an arbitrarily small negative number, in which case the lower bound becomes meaningless. For such cases, we can impose the trivial lower bound $\log_2(N)$.

In the remainder of this section we sketch the proof of Theorem \ref{warmup2_thm}. The main focus of this paper is on the adaptive setting, but for the sake of completeness we also include the non-adaptive version.

%%%%%%%%%%%%%%%%%%%%%%%%%%%%%%%%%
\subsection{Proof of the $\mathrm{QC}$  Upper Bound of Theorem \ref{warmup2_thm}} 
\label{QC_upper}
\begin{proof}
Direct application of Theorem \ref{warmup1_thm}.
\end{proof}
%%%%%%%%%%%%%%%%%%%%%%%%%%%%%%%%%
\subsection{Proof of the $\mathrm{QC}$  Lower Bound of Theorem \ref{warmup2_thm}} 
\label{QC_lower}
\begin{proof}
 Here we only consider the $1-\gamma_{qc}=\Theta(1)$ case when $\eta=1$. By Remark \ref{rem1}, the $1-\gamma_{qc}=o(1)$ case will follow from the $\mathrm{SQC}$ lower bound. Let $r\le(1-\epsilon)\frac{\log(N)}{\log(1/\gamma_{qc})}$ and $\epsilon \gg \frac{\log\log(N)}{\log(N)}$ with $\epsilon \rightarrow 0$. Let $Y$ be the number of subsets $W \subset \mathcal{R}$ with $|W|\le r$ for which $A_{W,\mathcal{C}}$ has no repeated columns. For the lower bound to hold we must show that $Y=0$ a.a.s.

Let us now select a set $R \subset \mathcal{R}$ of $r$ rows in advance, and when $A$ is revealed, let $\mathcal{A}_R$ be the event that $A_{R,\mathcal{C}}$ has no repeated columns. Then
\begin{equation}
\label{QC_eq1}
\P(Y>0) \le \E[Y] \le N^r\P(A \in \mathcal{A}_R).
\end{equation}
Using our result in equations \eqref{suen} and \eqref{lamdel} and because $1-\gamma_{qc}=\Theta(1)$ implies $\epsilon \ll q(1-\epsilon)$, for $N$ large enough

\begin{equation}
\label{QC_eq2}
\P(A \in \mathcal{A}_R)\le \exp(-\lambda + \Delta e^{2\delta}) < \exp\left(N^{2\epsilon}\left(- \frac14+2N^{\epsilon - \frac{3}{2}q(1-\epsilon)} \right)\right)< \exp\left(-\frac{1}{8}N^{2\epsilon} \right).
\end{equation}
Then,
\begin{equation}
\label{QC_eq3}
\E[Y] \le N^r \exp\left(-\frac{1}{8}N^{2\epsilon}  \right)\le \exp\left((1-\epsilon)\frac{\log^2(N)}{\log(1/\gamma_{qc})}-\frac18N^{2\epsilon}\right) \rightarrow 0
\end{equation}

since by assumption $\frac{1-\epsilon}{\log(1/\gamma_{qc})}=\Theta(1)$ and since $\epsilon \gg \frac{\log\log(N)}{\log(N)}$ implies $\log^2(N) \ll N^{2\epsilon}$. Finally, by equations \eqref{QC_eq1}  and \eqref{QC_eq3}, we have $Y=0$ a.a.s.
\end{proof}

%%%%%%%%%%%%%%%%%%%%%%%%%%%%%%%%%
\subsection{Proof of the $\mathrm{SQC}$  Upper Bound of Theorem \ref{warmup2_thm}} 
\label{SQC_upper}
In order to prove this upper bound, we analyse the performance of a greedy query selection algorithm called MAX-GAIN.

\begin{defn}[$k$-reducer]
\label{reducer} For a query $w \in \mathcal{R}$ and an observation $l \in \{0, 1\}$, let the targets agreeing with the pair $(w,l)$ be denoted as 
\begin{equation}
\mathcal{S}_A(w,l)=\{ v \in \mathcal{C} \mid A_{w,v}=l\}.
\end{equation}
Given an integer $k$ and the triple $(A,R_{j},A_{R_{j},\{v^\star\}})$, a row $w$ is called a \textit{$k$-reducer} if after adding $w$ to $R_{j}$, the worst case cardinality of $R_{j+1}$ is upper bounded by $k$, that is
\begin{equation}
\max\limits_{l \in \{0,1\}}|\mathcal{T}_{R_{j}}\cap \mathcal{S}_A(w,l) |  \le k.
\end{equation}
\end{defn}

\begin{defn}[MAX-GAIN]
\label{MAX-GAIN}
The \textit{MAX-GAIN} algorithm finds the target by iteratively selecting as a query the $k$-reducer with the smallest $k$. That is,
$$\mathrm{MAXGAIN}(A,R_{j},A_{R_{j},\{v^\star\}}) = \argmin\limits_{w\in V\setminus R_{j}}\max\limits_{l \in \{0,1\}} |\mathcal{T}_{R_{j}}\cap \mathcal{S}_A(w,l) |.$$
\end{defn}

Note that if $\frac{\log(N)}{N} \ll 1-\gamma_{sqc}$, there are a.a.s. no two identical columns in $A$ by Remark \ref{thm1_implies_p3}, in which case the MAX-GAIN algorithm always finds the target in at most $N$ steps. Moreover, if we can always find better reducers, the number of steps decreases dramatically. Since each node is connected to a $\gamma_{sqc}>1/2$ fraction of the nodes, it is reasonable to expect that we can find $k$-reducers with $k\approx |\mathcal{T}_{R_j}|\gamma_{sqc}$. The existence of such reducers would already imply the result we need.

\begin{lem}
\label{prog}
If MAX-GAIN can select a $(|\mathcal{T}_{R_j}|\gamma_{sqc}+f(|\mathcal{T}_{R_j}|))$-reducer in the $(j+1)^{th}$ step of the algorithm with $f(n)=o\left(\frac{n}{\log(n)}\right)$ for any $j\in \mathbb{N}$ for which the candidate set size is $|\mathcal{T}_{R_j}| = \Omega \left(\frac{\log(N)}{\log{(1/\gamma_{sqc})}}\right)$, then the algorithm finds the source in $(1+o(1))\frac{\log(N)}{\log{(1/\gamma_{sqc})}}$ steps.
\end{lem}

The proof of Lemma \ref{prog} is included in Section \ref{proof_l2}. For the $\mathrm{SQC}$ upper bound we will be able to prove the existence of a $(|\mathcal{T}_{R_j}|\gamma_{sqc}+f(|\mathcal{T}_{R_j}|))$-reducer for any candidate size, hence this condition of Lemma \ref{prog} may be ignored for the moment. The condition on the minimum candidate set size for which there exists a  $(|\mathcal{T}_{R_j}|\gamma_{sqc}+f(|\mathcal{T}_{R_j}|))$-reducer will be important in Section \ref{p1ub}.

Now we need to show the existence of such reducers. This will be a structural result on the matrix $A$ which holds independently of the state of our algorithm, so we find it useful to define another notion quite similar to $k$-reducers.

\begin{defn}[$f$-separator]
Let $f(n)\in \mathbb{N}\rightarrow\mathbb{R}^+$ be a function, and $\gamma_{sqc}$ as defined in Theorem \ref{warmup2_thm}. A set of columns $W\subseteq \mathcal{C}$, $|W|=n$ has an $f$\textit{-separator} if there is a row $w \in \mathcal{R}$ such that 

\begin{equation}
\max\limits_{l \in \{0,1\}}|W\cap \mathcal{S}_A(w,l) |\le n\gamma_{sqc}+f(n).
\end{equation}
\end{defn}

\begin{rem}
An $f$-separator $w$ for $\mathcal{T}_{R_j}$ is an $(|\mathcal{T}_{R_j}|\gamma_{sqc}+f(|\mathcal{T}_{R_j}|))$-reducer for the triple $(A,R_{j},A_{R_{j},\{v^\star\}})$. The difference between the two terms is that the term $f$-separator refers to a property of $A$ and $W$, whereas the term $k$-reducer refers to a property of the state of an algorithm. The role of these two terms is also quite different. A $k$-reducer with a small $k$ makes MAX-GAIN more efficient, whereas an $f$-separator with a small $f$ is a typical separator, and its existence makes the analysis of this upper bound easier. 
\end{rem}

\begin{proof}[Proof of the $\mathrm{SQC}$ upper bound of Theorem \ref{warmup2_thm}]
To use Lemma \ref{prog}, we have to show that for every $W \subseteq \mathcal{C}$ we have an $f$-separator with $f(n)=o\left(\frac{n}{\log(n)}\right)$.  Let $X_w=|W\cap \mathcal{S}_A(w,1)|$ and note that $\E[X_w]=nq$. It is clear that if $X_w$ is close to its expected value then $v$ must be an $f$ separator. Indeed, $|X_w-nq| \le f(n)$ implies
\begin{equation}
\label{abs_sep1}
X_w-nq \le f(n) \Rightarrow X_w \le nq+ f(n) \le n\gamma_{sqc} +f(n)
\end{equation}
and
\begin{equation}
\label{abs_sep2}
-(X_w-nq)-n+n \le f(n)\Rightarrow n-X_w \le n(1-q)+ f(n) \le n\gamma_{sqc} +f(n),
\end{equation}
hence $w$ is an $f$ separator. We first show that for any $v\in \mathcal{R}$ we have $|X_w-nq| \le f(n)$ with constant probability. Using Lemma \ref{chernoff} and substituting $f(n)=\sqrt{3n}$, we get
\begin{align}
\label{e-1}
\P( |X_w-nq| \ge f(n) ) &= \P\left( |X_w-\E[X_w]| \ge \frac{f(n)}{nq}nq \right) \le 2e^{\frac{-2nq\frac{f(n)^2}{n^2q^2}}{3}}=2e^{-\frac{6}{3q}}<e^{-1},
\end{align}
because $q\le1$. Since the random variables $X_w$ are mutually independent, the probability that none of the rows are $f$-separators for $W$ is upper bounded by $e^{-N}$. Let $Y$ be the number of subsets $W$ that do not have a $\sqrt{3n}$-separator. Then,
\begin{align*}
\E[Y] < \sum\limits_{W \subseteq \mathcal{C} } e^{-N} \le 2^Ne^{-N}   \rightarrow 0.
\end{align*}

By the first moment method we can conclude that every $W \subseteq \mathcal{C}$ has a $\sqrt{3n}$-separator a.a.s. By Lemma \ref{prog} this concludes the proof. 
\end{proof}

%%%%%%%%%%%%%%%%%%%%%%%%%%%%%%%%%
\subsection{Proof of the $\mathrm{SQC}$  Lower Bound of Theorem \ref{warmup2_thm}} 
\label{SQC_lower}
For this lower bound, we look for columns with identical elements similarly to Case 2 of the Proof of Theorem \ref{warmup1_thm}, part (ii). We have seen in Remark \ref{thm1_implies_p3}, that if $\frac{\log(N)}{N} \ll 1-\gamma_{sqc}$ then $A$ a.a.s. does not have two identical columns. However, any $r \times N$ of its submatrices will have two rows with identically 0 or 1 elements if $r\le (\eta-\epsilon) \frac{\log(N)}{\log(1/\gamma_{sqc})}$ with $\epsilon \gg \frac{\log\log(N)}{\log(N)}$, no matter how we select the rows (queries). Of course, it may differ which are the two columns that are identical based on our query selection. Recall, that at the end of Section \ref{definitions} we modelled the SMD as the number of steps in a two-player game, if both players play optimally. In this proof, we are essentially analysing a strategy for Player 2, who does not decide the target in advance, and always provides observation is $0$ if $q\le \frac12$ and $1$ if $q> \frac12$. By showing that with high probability any $r \times N$ submatrix of $A$ has at least two columns with identically 0 or 1 elements, we assure that Player 2 can follow this simple strategy, and the size of the candidate set will be at least two after $r$ queries, independently of the strategy of Player 1.

\begin{proof}[Proof of the $\mathrm{SQC}$  lower bound of Theorem \ref{warmup2_thm}]
Similarly to Section \ref{QC_lower}, let $Y$ be the number of subsets $W \subset \mathcal{R}$ with $|W|\le r$ for which $A_{W,\mathcal{C}}$ has at most one column with identically $0$ (if $q\le\frac12$) or $1$ (if $q>\frac12$) elements. For the lower bound to hold we must show $Y=0$ a.a.s. Let us now select a $R \subset \mathcal{R}$ of size $r$ in advance, and when $A$ is revealed let $\mathcal{A}_R$ be the event that $A_{R,\mathcal{C}}$ has at most one column with identical elements. Then,

\begin{equation}
\label{SQC_end}
\P(Y>0) \le E[Y] \le N^r\P(A \in \mathcal{A}_R).
\end{equation}

Let $Z_R$ be the number of identically $0$  (or $1$ if $q>\frac12$) columns in $A_{R,\mathcal{C}}$. By equation~\eqref{EZ},

\begin{equation}
\P(A \in \mathcal{A}_R) \le 2e^{\frac43}e^{-\frac{\E[Z_R]}{3}}
\end{equation}

Then, using the equation $\E[Z_R]=N\gamma_{sqc}^r$ and the definitions of $r$ and $\eta$,

\begin{align}
\label{SQC_eq3}
\E[Y]& \le N^r 2e^{\frac43}e^{-\frac{\E[Z_R]}{3}}\nonumber \\
&=2e^{\frac43}\exp\left(r\log(N)-\frac13N\gamma_{sqc}^r\right)\nonumber \\
&\le 2e^{\frac43}\exp\left(r\log(N)-\frac13N\gamma_{sqc}^{(\eta-\epsilon) \frac{\log(N)}{\log(1/\gamma_{sqc})}}\right)\nonumber \\
&= 2e^{\frac43}\exp\left(r\log(N)-\frac13N^{1-(\eta-\epsilon) }\right)\nonumber \\
&\le 2e^{\frac43}\exp\left(\frac{(\eta-\epsilon)\log^2(N)}{\log(1/\gamma_{sqc})}-\frac13N^{1-(1+\log_N\log(1/\gamma_{sqc}))+\epsilon}\right)\nonumber \\
&\le 2e^{\frac43}\exp\left(\frac{\log^2(N)-\frac13N^{\epsilon}}{\log(1/\gamma_{sqc})}\right) \rightarrow 0
\end{align}
since $\frac{1}{\log(1/\gamma_{sqc})}>1$ and  $\log^2(N)-\frac13N^{\epsilon}\rightarrow -\infty$ as long as $\epsilon \gg \frac{\log\log(N)}{\log(N)}$.  Finally, by equations \eqref{SQC_end}  and \eqref{SQC_eq3} we have $Y=0$ a.a.s.
\end{proof}
         
\section{Expansion Properties of $\mathcal{G}(N,p)$}
\label{metric}

\begin{figure}[h]
\begin{center}
  \includegraphics[width=0.7\textwidth]{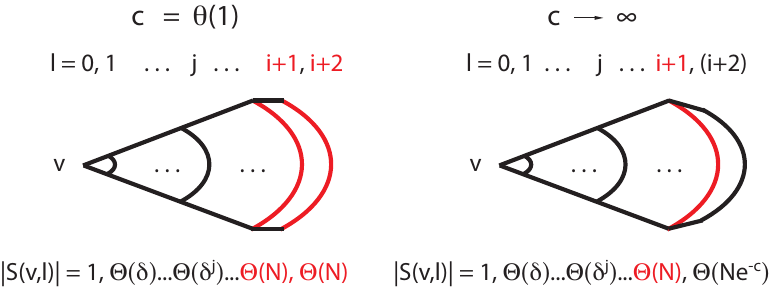}
  \caption{The arcs in the figures represent the level sets $\mathcal{S}_G(v,l)$ of $\mathcal{G}(N,p)$ for different ranges of $c$. The layers containing a constant fraction of the nodes marked red. In the $c=\Theta(1)$ case there are two such layers, whereas in the $c\rightarrow \infty$ there is only one. In this latter case the $(i+2)^{th}$ layer is in parenthesis because that layer may or may not exist depending on $c$. }
  \label{level_sets}

  \end{center}
\end{figure}

Before we proceed to our main results, we must establish some properties about the exponential growth of Erd\H{o}s-R\'enyi graphs in the sizes of the level sets $\mathcal{S}_G(v,l)$ defined below. This exponential growth is depicted on Figure \ref{level_sets}. Most statements in this section already appeared in \cite{bollobas2012metric} with a different notation, or can be easily derived from their results.

\begin{table*}[h!]
\begin{center}
{\renewcommand{\arraystretch}{1}%
\small{
\begin{tabular}{ |c|c|c|c|c|c| }
  \hline  $c(N)$ & \multicolumn{2}{c|}{$c=\Theta(1)$} & \thead{$1\ll c$ \\ $ c\ll  \log(\frac{N}{\delta^i})$} & \thead{ $\log(\frac{N}{\delta^i})\ll c$ \\ $c\ll 2\log(N)$ } & \thead{$2\log(N) \ll c$ \\$c \ll \delta$} \\ \hline \hline 
   $D$ & $i+2$ & $i+2$ & $i+2$ & $i+2$ & $i+1$ \\ \hline
   $|\mathcal{S}_G(v,i-1)|$ & $\delta^{i-1}$ &  $\delta^{i-1}$ & $\delta^{i-1}$ & \cellcolor{gray!25} $\delta^{i-1}$ & \cellcolor{gray!25} $\delta^{i-1}$\\ \hline
   $|\mathcal{S}_G(v,i)|$ & \cellcolor{gray!25} $\delta^i$ & \cellcolor{gray!25} $\delta^i$ &\cellcolor{gray!25} $\delta^i$ & \cellcolor{pink!25} $\delta^i$ & \cellcolor{pink!25} $\delta^i$ \\ \hline
   $|\mathcal{S}_G(v,i+1)|$ & \cellcolor{pink!25} $(1-e^{-c})N$ & \cellcolor{red!25} $(1-e^{-c})N$ & $\cellcolor{red!25} (1-e^{-c})N$ & $\cellcolor{red!25} \left(1-\frac{\delta^i}{N}\right)N$ &\cellcolor{red!25} $\left(1-\frac{\delta^i}{N}\right)N$ \\ \hline
   $|\mathcal{S}_G(v,i+2)|$ & \cellcolor{red!25} $e^{-c}N$ & \cellcolor{pink!25} $e^{-c}N$ & \cellcolor{pink!25} $e^{-c}N$ & $e^{-c}N$ & $0$ \\ \hline \hline

   MD ub \cite{bollobas2012metric} & \multicolumn{2}{c|}{T3.1 case 1} & T3.1 case 2.1 &  \multicolumn{2}{c|}{T3.1 case 2.2} \\ \hline
   MD lb \cite{bollobas2012metric} & \multicolumn{2}{c|}{ T4.2} & \thead{T4.2 ,\\ T4.4 case 1} & T4.4 case 2 & T4.3 \\ \hline
%   Our MD ub & \multicolumn{4}{c|}{Reduction to SQC ub} \\ \hline
 %  Our MD lb & Reduction to SQC lb  & \multicolumn{3}{c|}{Reduction to DQC lb} \\ \hline
    SMD ub & \multicolumn{2}{c|}{Similar to the SQC ub} & \multicolumn{3}{c|}{Use MD ub}  \\ \hline
    SMD lb & \multicolumn{5}{c|}{Coupling and an analysis similar to the DQC lb}  \\ \hline

\end{tabular}}}
\end{center}
  \caption{Overview of the main tools to prove Theorem \ref{main_thm}. Each column corresponds to a different range of parameter~$c$. The $c=\Theta(1)$ columns are split in two sub-columns: in the first, $e^{-c}>1-e^{-c}$ and in the second, $e^{-c}<1-e^{-c}$. Only the leading terms of the size of the level sets $S$ are shown. The largest level set is colored in red, and the second largest is colored in pink. The last level set before one of the two dominating level sets is colored in grey.  The bottom half of the table points to the proof of the upper/lower bound for each parameter range of Theorem \ref{main_thm}, both in previous work and in this paper.}
  \label{tools}
\end{table*}

\begin{defn}[level sets]
\label{levelsetdef}
For a graph $G=(V,E)$ and a node $v \in V$, let the \textit{level set} of $v$ be defined as $\mathcal{S}_G(v,l)=\{ w \in V \mid d(v,w)=l\}$ for every $l \in \{0, \dots, |V|\}$. The level sets from a set of nodes $V' \subseteq V$ is defined as $\mathcal{S}_G(V',l)=\bigcup_{v\in V'} \mathcal{S}_G(v,l).$
\end{defn}

We also define three functions $\delta, i$ and $c$ (all depending on the parameters of an Erd\H{o}s-R\'enyi graph $\mathcal{G}(N,p)$), which will be useful throughout the rest of this paper.

\begin{defn}[parameters $\delta, i$ and $c$ of the expansion properties]
\label{parameters}
Let $\delta=Np$, let $i\ge0$ be the largest integer such that $\delta^i=o(N)$, and finally let $c=\frac{\delta^{i+1}}{N}=\delta^ip$.
\end{defn}

In this paper we only consider connected Erd\H{o}s-R\'enyi graphs with $\delta=Np \gg \log(N)$. We defer the interpretation of these definitions and introduce the main technical lemma that establishes the exponential growth of the level sets. This lemma also appeared in Lemma 2.1 of \cite{bollobas2012metric} (we replaced their $O(1/\sqrt{\omega}) + O(d^i/n)$ term with $O(\zeta)$ for simplicity), with an extra condition which we removed (proof in Section \ref{proofs}).

\begin{lem}[Expansion property]
\label{expansion}With parameters $i,c$ and $\delta\gg \log(N)$ as defined in Definition \ref{parameters}, let $\zeta = \zeta(N)$ be a function tending to slowly to zero with $N$ such that 
\begin{equation}
\label{zeta1}
\zeta \ge \max\left( \sqrt{\frac{\log(N)}{\delta}}, \frac{\delta^i}{N}\right).
\end{equation}
For a node $v\in V$, let $\mathcal{E}(v,j)$ be the event that for every $l\le j$
\begin{equation} 
\label{expansioneq}
|\mathcal{S}_G\left(v, l\right)| = \left(1 + O\left(\zeta \right)\right)\delta^l,
\end{equation}
and for two nodes $v\ne u \in V$, let $\mathcal{E}_2(u,v,j)$ be the event that for every $l\le j$
\begin{equation} 
\label{expansioneq2}
|\mathcal{S}_G\left(\{u,v\}, l\right)| = 2\left(1 + O\left(\zeta \right)\right)\delta^l.
\end{equation}

For a subset $V' \subset V$ let $\mathcal{E}(V',j)=\bigcap_{v \in V'} \mathcal{E}(v,j)$ be the event that expansion properties hold for all nodes in $V'$, and let $\mathcal{E}(j)$ be a shorthand for $\mathcal{E}(V,j)$. Similarly, let $\mathcal{E}_2(j)= \bigcap_{u \ne v \in V} \mathcal{E}_2(u,v,j)$. Then, for $G$ sampled from $\mathcal{G}(N, p)$ the event $\mathcal{E}(i) \cap \mathcal{E}_2(i) $ holds a.a.s.
\end{lem}

\begin{cor}
\label{cor_S}
For every $v \in V$ we have $\sum_{l=1}^i |\mathcal{S}_G(v,l)|= (1+O(\zeta))\delta^i=o(N)$ a.a.s.
\end{cor}

\begin{proof}
By Lemma \ref{expansion} by taking $l\le i$ since $\zeta \gg 1/\delta$ and
\begin{align}
&\sum\limits_{l=1}^i |\mathcal{S}_G(v,l)| = \sum\limits_{l=1}^i \left(1 + O(\zeta) \right)\delta^l = (1+O(\zeta))\delta^i = o(N).
\end{align}
\end{proof}

Parameter $\delta$ is essentially the expected degree of each node in $G \sim \mathcal{G}(N,p)$. We require $\delta \ge \log(N)/\zeta^2 \gg \log(N)$, so the graph is a.a.s. connected. The function $\zeta$ serves as the error term. Note that equation \eqref{zeta1} implies that $\zeta\ge p$ for $i>0$. Parameters $i$ and $c$ are both derived from $1/\log_N(\delta)$; parameter $c$ is $\delta$ raised to one minus the fractional part of $\delta$, and parameter $i$ is the integer part of $1/\log_N(\delta)$, or more precisely the ceiling minus one. Qualitatively, the level set structure of $\mathcal{G}(N, p)$ has a periodic behaviour as we tune $p$. As $p$ decreases, in each such ``period'' the outmost layer gains more and more nodes until it is fully saturated and another layer appears. Roughly speaking, parameter $i$ indicates the ``period'', and parameter $c$ provides a fine-grain tuning of $p$ within a ``period''. However, the ``periods'' and the appearance of new level sets are not exactly aligned. The next lemma tells us about how the diameter depends on $\delta$ and $N$.

%%%%%%%%%%%%%%%%%% DIAMETER
\begin{lem}[Lemma 4.1 \cite{bollobas2012metric}]
\label{diameter}
Suppose that $\delta = pN \gg \log(N)$, and that for some integer $D\ge1$
\begin{equation}
\label{diameqs}
\frac{\delta^{D-1}}{N} - 2 \log(N) \rightarrow -\infty \qquad \text{and} \qquad \frac{\delta^D}{N} - 2 \log(N) \rightarrow \infty
\end{equation}
Then the diameter of $G$ sampled from $\mathcal{G}(N, p)$ is equal to $D$ a.a.s.
\end{lem}

\begin{cor}
\label{cor_diam}
Let $\mathcal{D}$ be the event that the diameter of $G$ sampled from $\mathcal{G}(N,p)$ is either $i+1$ or $i+2$ with all parameters, including $i$, given by Definition \ref{parameters}, and as always $\delta \gg \log(N)$. Then $\mathcal{D}$ holds a.a.s. \end{cor}

\begin{proof}
We distinguish three cases:
\begin{enumerate}
\item If $\delta^{i}/N - 2 \log(N) \rightarrow -\infty$ and $\delta^{i+1}/N - 2 \log(N) \rightarrow \infty$,
then taking $D=i+1$ in Lemma \ref{diameter} implies that the diameter is $i+1$ a.a.s.
\item If $\delta^{i+1}/N - 2 \log(N) \rightarrow -\infty$ and $\delta^{i+2}/N - 2 \log(N) \rightarrow \infty$,
then taking $D=i+2$ in Lemma \ref{diameter} implies that the diameter is $i+2$ a.a.s.
\item If $\delta^{i+1}/N-2\log(N)=\Theta(1)$, then let us consider $G_1=\mathcal{G}(N,p\omega)$ and $G_2=\mathcal{G}(N,p/\omega)$ with $\omega \rightarrow \infty$ very slowly. Using Lemma \ref{diameter}, for $\omega$ growing slowly enough, the graphs $G_1$ and $G_2$ have diameter $i+1$ and $i+2$ respectively a.a.s. Since having diameter at least $D$ is a monotone graph property, $G$ must also have diameter $i+1$ or $i+2$ a.a.s.
\end{enumerate}

There are no other cases than the three outlined above because the equations
\begin{align}
&\frac{\delta^{i}}{N} - 2 \log(N) \rightarrow -\infty 
\end{align}
and
\begin{align}
&\frac{\delta^{i+2}}{N} - 2 \log(N)= c\delta -2\log(N)\ge \left(\frac{c}{\zeta^2}-2\right)\log(N)\rightarrow \infty.
\end{align}
must always hold by Definition \ref{parameters} and the assumption $\delta \gg \log(N)$.
\end{proof}

The previous results shows that most of the nodes are either distance $i+1$ or $i+2$ away from any arbitrary node $v \in V$. We now extend Lemma \ref{expansion} to these level sets too.

\begin{lem}
\label{fraction}
For every $v \in V$, let $\mathcal{E}(v,i+1)$ be the intersection of $\mathcal{E}(v,i)$ and of the event
\begin{equation}
\label{aasSip1}
|\mathcal{S}_G(v,i+1)| =
\begin{cases} \left(1+O\left(\sqrt{\frac{\log(N)}{N}}\right)\right)Np &\mbox{if } i = 0 \\ 
\left(1-\left(e^{-c}+\frac{\delta^i}{N}\right)+O\left(\zeta\left(e^{-c}+\frac{\delta^i}{N}\right)+\sqrt{\frac{\log^2(N)}{N}}\right)\right) N& \mbox{if } i > 0  .
\end{cases}
\end{equation}

For a subset $V' \subset V$ let $\mathcal{E}(V',i+1)=\bigcap_{v \in V'} \mathcal{E}(v,i+1)$ be the event that expansion properties hold for all nodes in $V'$, and let $\mathcal{E}(i+1)$ be a shorthand for $\mathcal{E}(V,i+1)$. Then event $\mathcal{E}(i+1)$ holds a.a.s.

\end{lem}

\begin{proof}
For a fixed node $v\in V$, let us expose all of its edges (i.e., sample the edges of $\mathcal{G}(N,p)$ adjacent to $v$ in any order and reveal them), and do the same for each of its neighbors recursively until depth $i$. This way of exposing edges also exposes all nodes in $W=\bigcup_{l\le i}\mathcal{S}_G(v,l)$. 

After exposing these edges, we have for both $i=0$ and $i>0$ that
\begin{equation}
\label{SGip1}
|\mathcal{S}_G(v,i+1)|= \mathrm{Binom}\left( |V \setminus W|, 1-(1-p)^{|\mathcal{S}_G(v,i)|} \right),
\end{equation}
because $\mathcal{S}_G(v,i+1) =\{ w\in V \setminus W \mid \exists v' \in \mathcal{S}_G(v,i) \text{ s.t. }  d(v',w)=1\}$.

(i) In the $i>0$ case, let us condition on the event $\mathcal{E}(\{v\},i)$. By Corollary \ref{cor_S}, the set $V \setminus W$ has $N-(1+O(\zeta))\delta^{i}$ nodes. Then, we have

\begin{align}
\label{condExp}
\E[|\mathcal{S}_G(v,i+1)| \mid \mathcal{E}(\{v\},i)]&=(N-(1+O(\zeta))\delta^{i})(1-(1-p)^{|\mathcal{S}_G(v,i)|}) \nonumber \\
&\stackrel{\eqref{expansioneq}}{=}   (N-(1+O(\zeta))\delta^{i})(1-e^{-(p+O(p^2))\delta^i(1+O(\zeta))}) \nonumber \\
&=(N-\delta^{i})(1-e^{-c}(1+O(\zeta))) +O(\delta^i\zeta)\nonumber \\
& = N\left(\left(1-\frac{\delta^i}{N}\right)(1-e^{-c}) +O\left(\zeta e^{-c}+\frac{\delta^i}{N}\zeta \right)\right) \nonumber \\
& = N\left(1-(1+O(\zeta))\left(e^{-c}+\frac{\delta^i}{N}\right)\right).
\end{align}

Let us denote $\mu= \E[|\mathcal{S}_G(v,i+1)| \mid \mathcal{E}(\{v\},i)]$. Then, by Lemma \ref{chernoff} with $\tau=\sqrt{6\log(N)/\mu}$ we have

\begin{equation}
\label{mutau}
\P( \left| |\mathcal{S}_G(v,i+1)| - \mu \right| >  \tau\mu  \mid \mathcal{E}(\{v\},i)] ) < e^{-\frac{6\log(N)}{3}} =\frac{1}{N^2}.
\end{equation}

Since equation \eqref{condExp} imples $N/\log(N) \ll \mu$, we have $\tau<\sqrt{6\log^2(N)/N}$. This, together with equation \eqref{mutau} implies that for any $v\in V$, with probability at least $1-\frac{1}{N^2}$, we have
\begin{align}
|\mathcal{S}_G(v,i+1)| &=\left(1+O\left(\sqrt{\frac{\log^2(N)}{N}}\right)\right) N\left(1-(1+O(\zeta))\left(e^{-c}+\frac{\delta^i}{N}\right)\right)\nonumber \\
&=\left(1-\left(e^{-c}+\frac{\delta^i}{N}\right)+O\left(\zeta\left(e^{-c}+\frac{\delta^i}{N}\right)+\sqrt{\frac{\log^2(N)}{N}}\right)\right) N.
\end{align}
The desired result is implied by a union bound.

(ii) For $i=0$ we have $\E[\mathcal{S}_G(v,i+1)]=(N-1)p$. In this case, by Lemma \ref{chernoff} with $\tau=\sqrt{6\log(N)/N}$ we get  that with probability at least $1-\frac{1}{N^2}$ we have
\begin{align}
|\mathcal{S}_G(v,1)| &=\left(1+O\left(\sqrt{\frac{\log(N)}{N}}\right)\right) Np.
\end{align}
The desired result is again implied by a union bound.

%Finally, by the law of total expectation, Lemma \ref{expansion} and since $|\mathcal{S}_G(v,i+1)|\le N$,

%\begin{align}
%\E[|\mathcal{S}_G(v,i+1)|]&=\E[|\mathcal{S}_G(v,i+1)| \mid \mathcal{E}(\{v\},i)]\P(\mathcal{E}(\{v\},i))+\E[|\mathcal{S}_G(v,i+1)| \mid \overline{\mathcal{E}(\{v\},i)}]\P(\overline{\mathcal{E}(\{v\},i)})\nonumber \\
%&= N(1-e^{-c}) +o(N).
%\end{align}
\end{proof}

In Lemma \ref{fraction} we were quite precise about the error terms. A much weaker formulation of the same idea can be useful for interpreting Theorem \ref{main_thm}.

\begin{cor}
\label{top2}
In the same setting as in Lemma \ref{expansion}, the expected fraction of nodes in the level set with the largest expected size (conditioning on the expansion properties $\mathcal{E}(\{v\},i)$) is
$$
\begin{cases}(1+o(1))p &\mbox{if } i = 0 \mbox{ and }(1+o(1))p \ge 1-p \\ 
1-(1+o(1))p &\mbox{if } i = 0 \mbox{ and } p < 1-p \\ 
(1+o(1))e^{-c} & \mbox{if } i > 0 \mbox{ and } e^{-c} \ge 1-\left(e^{-c}+\frac{\delta^i}{N}\right)   \\
1-(1+o(1))\left(e^{-c}+\frac{\delta^i}{N}\right)& \mbox{if } i > 0  \mbox{ and } e^{-c} < 1-\left(e^{-c}+\frac{\delta^i}{N}\right).
\end{cases}$$
\end{cor}

\begin{proof}
The case $i=0$ is obvious. The case $i>0$  is a corollary of equation \eqref{condExp}. Indeed, the level set with the largest expected size as $N$ tends to infinity is $\mathcal{S}_G(v,i+2)$ if and only if $e^{-c} \ge 1-\left(e^{-c}+\delta^i/N\right)$, in which case it contains a fraction
$$1-\left(1-(1+o(1))\left(e^{-c}+\frac{\delta^i}{N}\right)\right)-o(1)=(1+o(1))e^{-c}$$ of all nodes by Lemma \ref{diameter} and  equation \eqref{condExp}. Otherwise, the level set with the largest expected size as $N$ tends to infinity is $\mathcal{S}_G(v,i+1)$, and its expected size is computed in equation \eqref{condExp}.

\end{proof}

The results in this section are summarised in the first five rows of Table \ref{tools}.
      
\section{Main Results}
\label{main_result}

We are ready to state our main theorem.
\begin{thm}
\label{main_thm}
Let $N \in \mathbb{N}$ and $p \in [0,1]$ such that $\frac{\log^5(N)}{N}\ll p$ and $\frac{1}{\sqrt{N}} \ll 1-p$. With the parameters given in Definition \ref{parameters}, let

\begin{equation}
\label{gammadef}
\gamma_{smd}=
\begin{cases}\max(p,1-p) &\mbox{if } i = 0  \quad (\text{i.e., }p=\Theta(1)) \\ 
\max(e^{-c},1-e^{-c}-\frac{\delta^i}{N}) & \mbox{if } i > 0 \quad (\text{i.e., }p=o(1)) \end{cases}
\end{equation}

and let 
\begin{equation}
\label{eq:etadef}
\eta=1+\log_N(\log(1/\gamma_{smd})).
\end{equation}
Finally, let $G=(V,E)$ be a realization of a $\mathcal{G}(N,p)$ random graph. Then, the following assertion holds a.a.s.

%\begin{enumerate}
%\item If $c=\Theta(1)$, then
%\begin{equation}
%\mathrm{SMD}(G)= (1+o(1))\frac{\log(N)}{\log(1/\gamma_{smd})}
%\end{equation}
%\item If $c\rightarrow \infty$, then
\begin{align}
& \mathrm{SMD}(G) \ge (\eta+o(1)) \frac{\log(N)}{\log(1/\gamma_{smd})} \nonumber \\ 
& \mathrm{SMDP}(G)  \le (1+o(1))\frac{\log(N)}{\log(1/\gamma_{smd})}.
\end{align}
%\end{enumerate}
\end{thm}

\begin{rem}
In case of $c\rightarrow \infty$ we have 
\begin{equation}
\label{1-g1pg}
\frac{1}{\log(1/\gamma_{smd})}=(1+o(1))\frac{1}{1-\gamma_{smd}} = (1+o(1))\left( e^{-c} + \frac{\delta^i}{N} \right)^{-1}.
\end{equation}
\end{rem}

\begin{rem}
\label{rem:boll_result}
The results for the MD of $\mathcal{G}(N,p)$ are of the same form (see Theorem 1.1 of \cite{bollobas2012metric}), but for the MD, instead of $\gamma_{smd}$, we have

\begin{equation}
\label{gammamddef}
\gamma_{md}=
\begin{cases} \sqrt{p^2+(1-p)^2} &\mbox{if } i = 0 \\ 
\sqrt{(e^{-c})^2 + (1 - e^{-c} - \frac{\delta^i}{N})^2} & \mbox{if } i > 0. \end{cases}
\end{equation}

Again, in case of $c\rightarrow \infty$ we have
\begin{equation}
\frac{1}{\log(1/\gamma_{md})} = (1+o(1))\left(e^{-c}  + \frac{\delta^i}{N}\right)^{-1} =(1+o(1))\frac{1}{\log(1/\gamma_{smd})}.
\end{equation}
\end{rem}

\begin{rem}
\label{rem:Fgammaeta}
The terms $F_\gamma$ and $F_\eta$ from Theorem \ref{meta_thm} can now be expressed based on Theorem \ref{main_thm} and Remark \ref{rem:boll_result}. With the definitions in equations \eqref{gammadef}, \eqref{eq:etadef} and \eqref{gammamddef}, taking the ratio of the appropriate lower and upper bounds, we can write

\begin{align}
\frac{\mathrm{SMD}(G)}{\mathrm{MD}(G)} &= (1+o(1))\frac{\log(1/\gamma_{md})}{\log(1/\gamma_{smd})} \quad \text{if } i=0\\
\frac{\mathrm{SMD}(G)}{\mathrm{MD}(G)} &\ge (\eta+o(1)) \quad \text{if } i>0,
\end{align}
hence we have
\begin{align}
F_\gamma &=  \frac{\log(\gamma_{md})}{\log(\gamma_{smd})}, \\ 
F_\eta &= \eta.
\end{align}
Finally, we note that since $p \in (0,1)$ and $e^{-c} \in (0,1)$, elementary analysis can show that $F_\gamma$ is a continuous function in $p$ taking every value in the interval $(1/2,1)$. A discussion on the value of $\eta$ is included in Section \ref{interpret}.
\end{rem}

%%%%%%%%%%%%%%%%%%%%%%%%%%%%%%%%%%%%%%%
\subsection{Connection between Theorems \ref{warmup2_thm} and \ref{main_thm}}
\label{interpret}

Clearly, there is a great deal of similarity between the MD/SMD in $\mathcal{G}(N,p)$ and the QC/SQC in $\mathrm{Ber}(q)^{N\times N}$. The final expression can always be written in the form $$(1+o(1))\frac{\log(n)}{\log(1/\gamma_{\star})},$$ where $\gamma_{\star}$ is a root mean square in the non-adaptive case and a maximum in the adaptive case. The main difference is that in binary search with randomly restricted queries, $\gamma_{\star}$ depends on parameter $q$ through a simple direct relation, whereas in random graph binary search the dependence of $\gamma_{\star}$ on $p$ is more complicated (see equation \eqref{gammadef}).

We can understand the mapping from $p$ to $\gamma_{\star}$, if we understand how we can map $p$ to the parameter $q$ in Theorem~\ref{warmup2_thm}, which we can do based on our results in Section \ref{metric}. When $p=\Theta(1)$, the mapping is just $p=q$. Indeed, since the diameter is 2 a.a.s, the vector $d(R,v)$ for $v\not\in R$ is essentially $\mathrm{Ber}(p)+1$. When $p=o(1)$, the size of either one or two level sets dominates the size of the others (see Figure $\ref{level_sets}$), hence the information we get is still basically a random Bernoulli vector, although here we must be more careful in the analysis. The mapping from $p$ to $q$ uses exactly the fraction of nodes in the largest level set established in Corollary \ref{top2}.

The $\eta$ term used in Theorem \ref{main_thm} serves the same purpose as in Theorem \ref{warmup2_thm}, the only difference is that in this case $\eta$ cannot be arbitrarily small (similarly to Remark 1.2 in \cite{bollobas2012metric}). Note that $i\ge1$ implies $\delta^i/N \ge (\log(N)\sqrt{N})^{-1}$, because otherwise $\delta^{2i}=o(N)$, which would imply that $i$ was not the largest integer for which $\delta^{i}=o(N)$, a contradiction with Definition \ref{parameters}. Since $-\log(1-x) \ge x$ and for $N$ large enough $1-\gamma_{smd} \ge e^{-c}+ \delta^i/N > \delta^i/N  \ge (\log(N)\sqrt{N})^{-1}$, we have 
\begin{align*}
\eta &= 1+\log_N(-\log(1-(1-\gamma_{smd})))\ge 1+\log_N(\gamma_{smd}))\\
&\ge 1-\log_N(\log(N)\sqrt{N})  = \frac12 -\frac{\log\log(N)}{\log(N)} \rightarrow \frac12.
\end{align*}
Since $\gamma_{smd}\ge 1/2$ implies $\log_N(\log(1/\gamma_{smd}))<0$, we have the bounds $1\ge \eta \ge 1/2$. Similarly to Section \ref{interpret_1_and_2}, we conjecture that the $1-\eta$ gap between the upper and lower bounds can be closed by improving the constant in front of the lower bound from $\eta$ to~$1$.

Since different ranges of parameters require different proof techniques both in this paper and in \cite{bollobas2012metric}, an overview of the proofs is presented in Table \ref{tools} for better clarity.

%%%%%%%%%%%%%%%%%%%%%%%%%%%%%%%%%%%%%%%%%
\subsection{Proof of Theorem \ref{main_thm} for $p=\Theta(1)$, Upper Bound}
\label{p1ub}
\begin{proof}
The $p=\Theta(1)$ case of Theorem \ref{main_thm} seems very similar to Theorem \ref{warmup2_thm} because entries of the distance matrix of $G \sim \mathcal{G}(N,p)$ are essentially $\mathrm{Ber}(p)+1$ random variables. However, the matrix is always symmetric, which causes some complications in the proof.

We start by providing an analogous definition of an $f$-separator for graphs.

\begin{defn}[$f$-separator]
\label{f-sep-graphs}
Let $f(n)\in \mathbb{N}\rightarrow\mathbb{R}^+$ be a function. A set of nodes $W\subseteq V$, $|W|=n$ has an $f$\textit{-separator} if there is a node $w \in V$ such that 
\begin{equation}
\max\limits_{l \in \mathbb{N}}|W\cap \mathcal{S}_{G}(w,l) |\le n\gamma_{smd}+f(n).
\end{equation}
\end{defn}

For the upper bound, the statement that for any set $W \subseteq V$, any node $w\in V$ can independently be an $f$-separator is not true anymore in contrast to the proof of Theorem \ref{warmup2_thm}, since the neighbourhoods of the nodes in $W$ are slightly correlated. However, the statement is still true for nodes $w \in V\setminus W$. Hence the proof will go on two steps. In step~1 we prove that $V$ has a $2\sqrt{n}$-separator a.a.s., and next in step 2 we prove that any set $W$ of cardinality at most $\gamma_{smd}N+2\sqrt{N}$ has an $f$-separator with $f(n)=o\left(\frac{n}{\log(n)}\right)$.

For step 1, let us pull aside from $V$ a random subset $F\subset V$ of cardinality $|F|=\log\log(N)$. By equation~\eqref{e-1} with $q=p$ and $X_w=|V' \cup \mathcal{S}_G(w,1)|$, each $w\in F$ is not a $\sqrt{3n}$-separator of $V'=V\setminus F$ with probability at most $e^{-1}$. Since these events are independent, the probability that no node $w\in F$ is an $\sqrt{3n}$-separator of $V'$ is then $e^{-\log\log(N)}=\log(N)^{-1}\rightarrow 0$. On the other hand, a $\sqrt{3n}$-separator of $V'$ is also a $2\sqrt{n}$-separator of $V$, since $\sqrt{3N}+\log\log(N)< 2\sqrt{N}$ for $N$ large enough.

For step 2, we repeat the calculation in equation \eqref{e-1} with $f(n)=\sqrt{\frac{6}{1-\gamma_{smd}}}$. Let $X_w=|W\cap \mathcal{S}_G(w,1)|$, then

\begin{align}
\label{e-12}
\P( |X_w-np| \ge f(n) )&= \P\left( |X_w-\E[X_w]| \ge \frac{f(n)}{np}np \right) \le 2e^{\frac{-2np\frac{f(n)^2}{n^2p^2}}{3}} \nonumber \\
&=2e^{-\frac{12}{3p(1-\gamma_{smd})}}<e^{-\frac{2}{1-\gamma_{smd}}},
\end{align}
because $p\le1$. Note that by equations \eqref{abs_sep1} and \eqref{abs_sep2} with $q=p$ and $\gamma_{sqc}=\gamma_{smd}$, the event $|X_w-np| \ge f(n)$ implies that $w$ is an $f$-separator for $W$.

Let $Y$ be the number of subsets $W$ of cardinality at most $\gamma_{smd}N+2\sqrt{N}$ that do not have a $\sqrt{\frac{6}{1-\gamma_{smd}}}$-separator. Then,

\begin{align*}
\P(Y>0) &\le \E[Y] < \sum\limits_{|W|\le \gamma_{smd}N+2\sqrt{N}} e^{-\frac{2}{1-\gamma_{smd}}(N-|W|)} \\
&\le 2^Ne^{-\frac{2}{1-\gamma_{smd}} (N-(\gamma_{smd}N+2\sqrt{N}))} \\
&< e^{-N+\frac{4\sqrt{N}}{1-\gamma_{smd}}}   \rightarrow 0,
\end{align*}

as long as $1-\gamma_{smd} \gg \frac{1}{\sqrt{N}}$. The existence of a $2\sqrt{n}$-separator for $V$ (step 1) and a  $\sqrt{\frac{6}{1-\gamma_{smd}}}$-separator for all $W\subseteq V$ of cardinality at most $\gamma_{smd}N+2\sqrt{N}$ (step 2) holds together a.a.s. by the union bound. Note that for $\gamma_{smd}\rightarrow 1$

\begin{align}
\label{longder}
\sqrt{\frac{6}{1-\gamma_{smd}}} &\stackrel{\eqref{1-g1pg}}{=}
\sqrt{\frac{6\left(1+o\left(1\right)\right)}{\log\left(\frac{1}{\gamma_{smd}}\right)}} \ll
\frac{1}{\log\left(\frac{1}{\gamma_{smd}}\right)\left(1+\log\left(\frac{1}{\log\left(\frac{1}{\gamma_{smd}}\right)}\right)\right)}\nonumber \\
& \ll \frac{\log\left(N\right)}{\log\left(\frac{1}{\gamma_{smd}}\right)\log\log(N)\left(\log\log\left(N\right)-\log\log\log(N)+\log\left(\frac{1}{\log\left(\frac{1}{\gamma_{smd}}\right)}\right)\right)} \nonumber \\
&= \frac{n}{\log\left(n\right)}
\end{align} 

for $n=\frac{\log(N)}{\log\log(N)\log(1/\gamma_{smd})}$. Hence, $\sqrt{\frac{6}{1-\gamma_{smd}}}=o\left(\frac{n}{\log(n)}\right)$ for $n = \Omega\left(\frac{\log(N)}{\log(1/\gamma_{smd})}\right)$, which is the necessary condition on $f(n)$ to apply Lemma \ref{prog}.

Since we were able to prove the existence of $(|\mathcal{T}_{R_j}|\gamma_{sqc}+f(|\mathcal{T}_{R_j}|))$-reducers for $j=1$ (step 1 of the analysis), and for any $j>1$ with candidate set size $|\mathcal{T}_{R_j}| = \Omega \left(\frac{\log(N)}{\log(1/\gamma_{smd})}\right)$ (step 2 of the analysis), Lemma \ref{prog} concludes the proof.  
\end{proof}
%%%%%%%%%%%%%%%%%%%%%%%%%%%%%%%%%%%%%%%%%
\subsection{Proof of Theorem \ref{main_thm} for $p=\Theta(1)$, Lower Bound}
\begin{proof}
The lower bound will be more straightforward given the SQC lower bound. We prove that Player 2 can follow the strategy of providing observation 0 if $p\le\frac12$ and 1 if $p >\frac12$, because every $r \times N$ submatrix of the distance matrix of $G$ has two columns, none of which has the same index as the indices of the rows, with entries identically equal to 1 or 2. We essentially sacrifice the entries that appear twice in the matrix (due to symmetry) for the ease of analysis. The entire derivation in the SQC lower bound follows except that now $\E[Z_R]=(N-r)\gamma_{smd}^r$, as there are only $N-r$ columns to chose from. Then, the rest of the computation is almost identical to equation \eqref{SQC_eq3}.

\begin{align}
\label{SMDp1lower}
P(Y>0)\le \E[Y]& \le N^r 2e^{\frac43}e^{-\frac{\E[Z_R]}{3}}\nonumber \\
&=2e^{\frac43}\exp\left(r\log(N)-\frac13(N-r)\gamma_{smd}^r\right)\nonumber \\
&\le 2e^{\frac43}\exp\left(r(\log(N)-\frac13(N-r)\gamma_{smd}^{(\eta-\epsilon) \frac{\log(N)}{\log(1/\gamma_{smd})}}\right)\nonumber \\
&= 2e^{\frac43}\exp\left(r(\log(N)+\frac13N^{-(\eta-\epsilon)})-\frac13N^{1-(\eta-\epsilon) }\right)\nonumber \\
&\le 2e^{\frac43}\exp\left(\frac{(\eta-\epsilon)\log(N)(\log(N)+1)}{\log(1/\gamma_{smd})}-\frac13N^{1-(1+\log_N\log(1/\gamma_{smd}))+\epsilon}\right)\nonumber \\
&\le 2e^{\frac43}\exp\left(\frac{\log^3(N)-\frac13N^{\epsilon}}{\log(1/\gamma_{smd})}\right) \rightarrow 0
\end{align}
since $\frac{1}{\log(1/\gamma_{smd})}>1$ and  $\log^3(N)-\frac13N^{\epsilon}\rightarrow -\infty$ as long as $\epsilon \gg \frac{\log\log(N)}{\log(N)}$.  

In the rest of this section we prove Theorem \ref{main_thm} for $p=o(1)$. Both the upper and lower bound will share some similarities with the $p=\Theta(1)$ case, but they will be much more involved. The reason that we still included the $p=\Theta(1)$ case is to study how far we can push the upper bound on $p$. In the end, our results on the SMD hold up to the bound $1-p\gg \frac{1}{\sqrt{N}}$. This is a slightly better bound then the  $1-p\ge \frac{3\log\log(N)}{\log(N)}$ up to which results of \cite{bollobas2012metric} on the~MD.
\end{proof}
\subsection{Proof of Theorem \ref{main_thm} for $p=o(1)$, Lower Bound}
%%%%%%%%%%%%%%%%%%%%%%%%%%%%%%%%
\begin{proof}
This proof is based on a coupling between the graph case and a simple stochastic process, which we can analyse similarly to the SQC lower bound. We start by upper bounding the probability that a random set is a resolving set by the probability that a certain survival process leaves at least two of the nodes alive. Since this survival process is still too complicated to analyse, we are going to introduce a second survival process later in the proof, which will give us the desired bound.

As usual, we first select queries $R=\{ w_1, \dots, w_r\}$ at random, with $|R| \le r=(\eta-\epsilon) \frac{\log(N)}{\log(1/\gamma_{sqc})}$ and $\epsilon$ slowly decaying to zero. In the proof of the SQC lower bound, we had an explicit lower bound on how slowly $\epsilon$ must tend to zero, however, this time we do not provide such guarantees for the sake of simplicity.

Let $l^\star$ be the index of the largest level set in expectation. Using the results of Corollary \ref{top2}, 
\begin{equation}
\label{jstar}
l^\star = \begin{cases} i+1 &\mbox{if } e^{-c} < 1-e^{-c}-\frac{\delta^i}{N} \\ 
i+2 & \mbox{if }  e^{-c} \ge 1-e^{-c}-\frac{\delta^i}{N} .\end{cases}
\end{equation}

Let $\mathcal{R}_R$ be the event that the randomly sampled set $R$ of size $|R|$ is a resolving set in $G$, and $\mathcal{R}$ that there exists at least one resolving set. Similarly to Section \ref{SQC_lower} we want to upper bound $\P(\mathcal{R})$ by $N^r \P(\mathcal{R}_R)$, and $\P(\mathcal{R}_R)$ by the probability of the event that there are at least two distinct nodes $u \ne v \in V$ with $d(R,u)=d(R,v)= l^\star \mathbf{1}$. 

Since $R$ is uniformly random, we may sample it before any of the edges in $G$ are exposed. Let us now expose the edges of $G$ similarly to the proof of Lemma \ref{fraction}, except this time starting from the set $R$ instead of a single node. Notice that before any of the graph is exposed, any of the nodes $v\in V\setminus R$ could possibly have $d(R,v)= l^\star \mathbf{1}$. Then, as more and more edges get exposed, many of the nodes lose this property. For instance the neighbors of the nodes in $R$ cannot have $d(R,v)= l^\star \mathbf{1}$, because $l^\star>i \ge1$. Hence focusing on the event $\mathcal{R}_R$, this exploration process of the graph can be seen as a survival process, where at least two nodes must survive.

\begin{defn}[ESP]
\label{ESP}
In the \textit{exploration survival process} (ESP) all nodes $v\in V\setminus R$ start out alive. In step $l<l^\star$, we expose all unexposed edges incident to the nodes in $\mathcal{S}_G(w_j, l)$ to expose $\mathcal{S}_G(w_j, l+1)$ for all $j \in \{1,\dots, r\}$. Every node exposed this way dies. Then, if $l^\star=i+1$ we play an extra round, in which we expose all unexposed edges incident to the nodes in $\mathcal{S}_G(w_j, i+1)$, and a node that was still alive at the end of round $l^\star-1$ survives this round if it connects to $\mathcal{S}_G(w_j, i)$ for all $j \in \{1,\dots, r\}$. If $l^\star=i+2$, there is no extra round.
\end{defn}

Note that the ESP always takes $i+1$ steps, it is only the nature of the final step that depends on $l^\star$. The event that $v$ survives the ESP is equivalent to $d(R,v)= l^\star \mathbf{1}$, unless $l^\star=i+2$ and event $\mathcal{D}$ in Corollary \ref{cor_diam} does not hold (in this case node $v$ surviving the ESP could have $d(w_j,v)> l^\star$). Since $\mathcal{D}$ holds a.a.s., we can assume it holds (formally, we may intersect all of our events with $\mathcal{D}$ and apply a union bound in the end).

In the first $l<l^\star$ rounds the probability of survival is $\rho^{(0)}_l=(1-p)^{|\mathcal{S}_G(R,l-1)|}$, and this probability is itself a random variable. When $l=l^\star=i+1$, we need each node to connect to each $\mathcal{S}_G(w_j,l-1)$, but these sets might have an intersection, so the exact value of the probability of survival (which we call $\rho^{(0)}_{l^\star}$) is complicated to write down. Fortunately, $\rho^{(1)}_{l^\star}$ can be lower bounded by $\prod_{j=1}^r \left( 1-(1-p)^{|\mathcal{S}_G(w_j,l^\star-1)|} \right)$, since the events of connecting to $\mathcal{S}_G(w_j,l-1)$ for each $j \in \{1,\dots, r\}$ are positively correlated. This motivates an alternative but still complicated survival process, which will serve as a bridge to the simple process we will finally analyse. 

\begin{defn}[CSP]
\label{CSP}
In the \textit{complex survival process} (CSP) all nodes $v\in V\setminus R$ start out alive. In each of the $i+1$ rounds, each node survives with probability $\rho^{(1)}_l$, where

\begin{equation}
\label{rho1}
\rho^{(1)}_l = \begin{cases} (1-p)^{|\mathcal{S}_G(R,l-1)|} & \mbox{if } l<l^\star \\
\prod\limits_{j=1}^r \left( 1-(1-p)^{|\mathcal{S}_G(w_j,l^\star-1)|} \right) & \mbox{if } l=l^\star \text{ and } l^\star=i+1,
\end{cases}
\end{equation}
where the sets $\mathcal{S}_G(R,l-1)$ are the same as in the ESP.
\end{defn}

Let $Y_0$ (respectively,  $Y_1$) be the indicator variable that at least two nodes survive the ESP (respectively,  CSP). Then $Y_0=0$ is the same event as $\mathcal{R}_R$ and $\rho^{(1)}_l \le \rho^{(0)}_l$, which implies $\P(\mathcal{R}_R)= \P(Y_0=0) \le \P(Y_1=0)$. However, the CSP is still too difficult to analyse even with the lower bound in \eqref{rho1} because each term is itself a random variable depending on earlier levels. Instead, we study a ``simple'' survival process.

\begin{cor}
\label{detcorollary}
Let us denote the event $\mathcal{E}(R,l^\star-1)$ by $\mathcal{E}_R$. Then, if $\mathcal{E}_R$ holds, there exists a constant $C$ such that
\begin{equation}
\label{detupper}
|\mathcal{S}_G(R,l-1)| \le r\delta^{l-1}\left(1+C\zeta\right)
\end{equation}
for all $R$ and $l<l^\star$, and when $l^\star=i+1$
\begin{equation}
\label{detlower}
|\mathcal{S}_G(w_j,l^\star-1)|  \ge \delta^{l^\star-1}\left(1-C\zeta\right)
\end{equation}
for all $w_j$.
\end{cor}

\begin{proof}
The existence of such a constant $C$ is implied by equation \eqref{expansioneq} in Lemma \ref{expansion} (for equation \eqref{detupper} since we need an upper bound, the intersection of the sets $\mathcal{S}_G(v,l-1)$ can be ignored).
\end{proof}

\begin{defn}[SSP]
\label{SSP}
In the \textit{simple survival process} (SSP) all nodes $v\in V\setminus R$ start out alive. In each of the $i+1$ rounds, each node survives with probability $\rho^{(2)}_l$, where
 \begin{equation}
\label{rho2}
\rho^{(2)}_l=\begin{cases} (1-p)^{r\delta^{l-1}\left(1+C\zeta\right)} & \mbox{if } l<l^\star \\
\left( 1-(1-p)^{(1-C\zeta)} \right)^r & \mbox{if } l=l^\star  \text{ and } l^\star=i+1,
\end{cases}
\end{equation}
and $C$ is the constant in Corollary \ref{detcorollary}.
\end{defn}

Let $Y_2$ be the indicator variable that at least two nodes survive the SSP. Equations \eqref{detupper} and \eqref{detlower} imply $\rho^{(1)}_l\ge\rho^{(2)}_l$, so it is in fact easier for nodes to survive the CSP than the SSP (the words ``simple'' and ``complex'' in the names of the terms SSP and CSP refer to the difficulty of analysis not the difficulty of survival). Hence, we should be able to prove $\P(Y_1=0, \mathcal{E}_R)\le \P(Y_2=0, \mathcal{E}_R)$, and we will prove it rigorously by coupling $(Y_1,\mathcal{E}_R)$ and $(Y_1,\mathcal{E}_R)$. Recall, that a coupling is a joint distribution $((\hat{Y}_1,\hat{\mathcal{E}}_R),(\hat{Y}_2,\hat{\mathcal{E}}_R))$ on $\{0,1\} \times \{0,1\} $ with the property that its first marginal is $(Y_1,\mathcal{E}_R)$ and the second marginal is $(Y_2,\mathcal{E}_R)$.

\begin{figure}[h!]
\begin{center}
  \includegraphics[width=\textwidth]{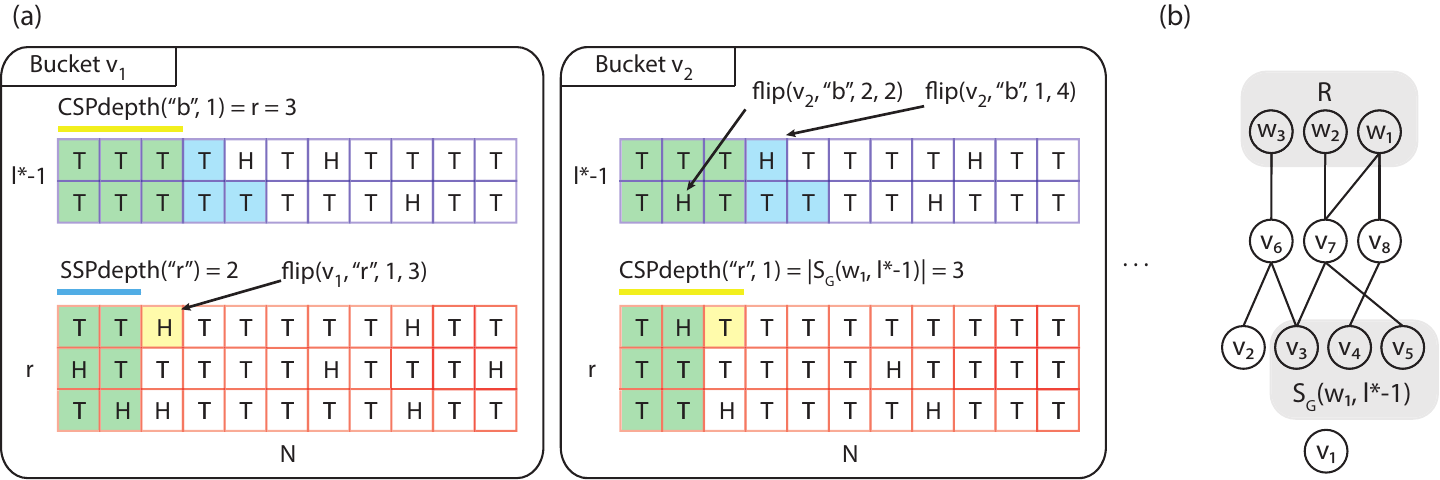}
  \caption{Part (a) of the figure shows a (partial) example sample of the elementary events of the coupled joint distribution.  Here we set $N=11, r=3$ and  $l^\star=i+1=3$. These parameters might not actually correspond to any pair of parameters $(N,p)$, we only use them for the example. The colors signify which survival processes use which coin flips; light blue is for the SSP, yellow is for the CSP and green is for the coin flips used by both of them. We should show $N-r=8$ buckets in total; one for each node $v_i \in V \setminus R$, but we only show the bucket for two nodes $v_1, v_2$ in the interest of space. Node $v_1$ survives the first $l^\star-1$ rounds in both processes, but survives the last round only in the CSP. In the SSP, it does not survive because the first red subbucket contains no head for this process (the CSP is saved by $\mathrm{flip}(v_1,"r",1,3)$). Node $v_2$ dies in the first round of the SSP (because of $\mathrm{flip}(v_2,"b",1,4)$) and in the second round of the CSP (because of $\mathrm{flip}(v_2,"b",2,2)$).\newline
Part (b) of the figure shows (a possible realisation of) the ESP corresponding the coin flips in part (a). The edges incident to nodes $v_1$ and $v_2$ correspond to the green coin flips in the blue subbuckets in part (a) of the figure. Only one such edge is present in this realization of the ESP; the edge between $v_2$ and $v_6$. This edge corresponds to $\mathrm{flip}(v_2,"b",2,2)$, which is indeed the only green head in the blue subbuckets in part (a) of the figure. Part (b) of the figure also explains the values of $\mathrm{CSPdepth}("r",j)$ in part (a). Indeed, we can check that $|\mathcal{S}_G(w_1,l^\star-1)|=3$ and $|\mathcal{S}_G(w_2,l^\star-1)|=|\mathcal{S}_G(w_3,l^\star-1)|=2$. Similarly, we can check that $|\cup_{j=1}^r \mathcal{S}_G(w_j,1)|= |\{ v_6, v_7, v_8 \} | =3$, which corresponds to the value of  $\mathrm{CSPdepth}("b",2)$ in part (a). Note that only coin flips in the blue subbuckets can correspond to edges in the ESP, as in the coupling we only simulate the ESP until round $l^\star-1$.}
  \label{buckets}
  \end{center}
\end{figure}

We define the joint distribution $(\hat{Y}_1,\hat{Y}_2,\hat{\mathcal{E}}_R)$ by specifying how to sample from it. We will simultaneously play the ESP, CSP and the SSP, and since all three processes can be simulated using Bernoulli trials with parameter $p$, we will use the same outcomes of these trials whenever possible. More precisely, the probability space of $(\hat{Y}_1,\hat{Y}_2,\hat{\mathcal{E}}_R)$ is the probability space of $N(N-r)(l^\star-1+r)$ coin flips where the probability of heads is $p$. The first $N(N-r)(l^\star-1+r)$ coin flips are first organized into $N-r$ buckets of size $N(l^\star-1+r)$ indexed by nodes $v \in V\setminus R$. Then, each of the $N-r$ buckets are further divided into $l^\star-1$ blue and $r$ red subbuckets, all of size $N$. The $k^{th}$ coin flip in the $l^{th}$ blue and respectively the $j^{th}$ red subbucket of the bucket corresponding to node $v$ is called $\mathrm{flip}(v,"b",l,k)$ and respectively $\mathrm{flip}(v,"r",j,k)$. Figure \ref{buckets} explains the structure and function of these buckets through an example.

To define $\hat{Y}_1$, we must explain how to simulate the CSP in Definition \ref{CSP} using the coin flips in the blue and red subbuckets. To simulate the CSP we must know the level set sizes $\mathcal{S}_G(R,l-1)$ for $l\le l^\star$, which requires simulating the ESP at least up to round $l^\star-1$. Note that in the ESP, we only expose edges between one dead and one alive node (if we consider the set $R$ dead at the start). For each such exposure of an edge between dead node $u$ and alive node $v$ in round $l<l^\star$, we use $\mathrm{flip}(v,"b",l,k)$, where $k$ is the lowest index for which $\mathrm{flip}(v,"b",l,k)$ has not been used so far in the ESP. This way, the exact mapping between edges and coin flips can change depending on the order we sample the edges in each round, but this will not affect the coupling. Until round $l^\star-1$, we are able to perfectly simulate the ESP, and fortunately we already have enough information to simulate the CSP. In round $l=l^\star$, we already know enough to finish simulating the CSP and we may ignore the ESP. We simply define $\hat{Y}_1$ to be the indicator of the event that there exists $V_{CSP}\subset V$ with $|V_{CSP}|\ge2$, such that for any $v \in V_{CSP}$ we have no head in the set
\begin{equation}
\{ \mathrm{flip}(v,"b",l,k) \mid 1 \le k  \le \mathrm{CSPdepth}("b",l) \} 
\end{equation}
for each positive integer $l < l^\star$, where
\begin{equation}
\mathrm{CSPdepth}("b",l) = \left|\bigcup\limits_{j=1}^r\mathcal{S}_G(w_j,l-1)\right| 
\end{equation}
This set is exactly the ``used'' nodes of the ESP, since the ESP and the CSP are identical in rounds $l<l^\star$. In addition, if $l^\star=i+1$, we also need that there is at least one head in the set
\begin{equation}
\{ \mathrm{flip}(v,"r",j,k) \mid 1 \le k  \le \mathrm{CSPdepth}("r",j)\} 
\end{equation}
for each positive integer $j \le r$, where
\begin{equation}
\mathrm{CSPdepth}("r",j) =  |\mathcal{S}_G(w_j,l^\star-1)|.
\end{equation}
It is clear that each bucket has at least as many coin flips as we need in each step, and that $\hat{Y}_1$ and $Y_1$ have the same distribution. Note that since we know the cardinality of the level sets $\mathcal{S}_G(R,l-1)$ for all $l\le l^\star$, we can also determine $\hat{\mathcal{E}}_R$.

The random variable $\hat{Y}_2$ can simply be defined as the indicator of the event that there exists $V_{SSP}\subset V$ with $|V_{SSP}|\ge2$, such that for any $v \in V_{SSP}$ we have no head in the set 
\begin{equation}
\{ \mathrm{flip}(v,"b",l,k) \mid 1 \le k  \le \mathrm{SSPdepth}("b",l) \}
\end{equation}
for each positive integer $l < l^\star$, where
\begin{equation}
\mathrm{SSPdepth}("b",l) =  r\delta^{l-1}\left(1+C\zeta\right).
\end{equation}
In addition, if $l^\star=i+1$, we also need that there is at least one head in the set
\begin{equation}
\{ \mathrm{flip}(v,"r",j,k) \mid 1 \le k  \le \mathrm{SSPdepth}("r")\} 
\end{equation}
for each positive integer $j \le r$, where
\begin{equation}
\mathrm{SSPdepth}("r") = (1-C\zeta)\delta^{l^\star-1}
\end{equation}
 (this depth does not depend on $j$, but it still needs to hold for $r$ subbuckets). Clearly, $\hat{Y}_2$ and $Y_2$ have the same distribution.

By equation \eqref{detupper}, if the event $\hat{\mathcal{E}}_R$ holds, each coin flip used by the CSP in rounds $\{1, \dots, i+1\}$ is also used by the SSP. Hence if a node survives in SSP it must also survive in the CSP. When $l^\star=i+1$, the situation is reversed in the $(l^\star)^{th}$ round. By equation \eqref{detlower}, each coin flip used by the SSP is used by the CSP. However, this time we need heads to survive, so again if a node survives in SSP it must also survive in CSP. Hence, $\hat{\P}(\hat{Y}_1,\hat{\mathcal{E}}_R)<\hat{\P}(\hat{Y}_2,\hat{\mathcal{E}}_R)$. We are now ready to use our coupling to bound the probability that there exists a resolving set. 

\begin{align}
\label{exres}
\P(\mathcal{R}) &\le \P(\mathcal{R} ,\mathcal{E}) + \P(\overline{\mathcal{E}}) \nonumber \\
& \le \sum\limits_{|R|\le r} \P(\mathcal{R}_R | \mathcal{E}) + \P(\overline{\mathcal{E}})  \nonumber \\
& \le \sum\limits_{|R|\le r}\P(Y_1=0 | \mathcal{E}) + \P(\overline{\mathcal{E}})  \nonumber \\
& \le N^r \frac{\P(Y_1=0,\mathcal{E}_{R})}{\P(\mathcal{E})} + \P(\overline{\mathcal{E}}) = N^r \frac{\hat{\P}(\hat{Y}_1=0,\mathcal{E}_{R})}{\P(\mathcal{E})}  + \P(\overline{\mathcal{E}}) \nonumber \\
& \le N^r \frac{\hat{\P}(\hat{Y}_2=0,\mathcal{E}_{R})}{\P(\mathcal{E})} + \P(\overline{\mathcal{E}}) \nonumber \\
& \le N^r \frac{\hat{\P}(\hat{Y}_2=0)}{\P(\mathcal{E})}  + \P(\overline{\mathcal{E}}) = N^r \frac{\P(Y_2=0)}{\P(\mathcal{E})} + \P(\overline{\mathcal{E}}).
\end{align}

Now we proceed to upper bounding $\P(Y_2=0)$. Let $Z_v$ be the indicator of the event that node $v \in V\setminus R$ survives in the SSP. We need to distinguish the two cases (i) $ e^{-c} \ge 1-e^{-c}-\frac{\delta^i}{N}$ and  (ii) $e^{-c} < 1-e^{-c}-\frac{\delta^i}{N}$. 

(i) In the case $e^{-c} < 1-e^{-c}-\frac{\delta^i}{N}$, since by equation \eqref{jstar} we have $l^\star=i+1$, 
\begin{align}
\label{caseip1}
\P(Z_v=1)& \ge \left(\prod\limits_{l=1}^{l^\star-1}(1-p)^{r\delta^{l-1}\left(1+C\zeta\right)}\right) \left( 1-(1-p)^{(1-C\zeta)\delta^{l^\star-1})} \right)^r \nonumber \\
&\ge e^{-p(1+o(1))r\delta^{i-1}}(1-e^{-p(1+o(1))\delta^{i}})^r \nonumber \\
&\ge \left(1-\frac{\delta^{i}}{N}\right)^{r(1+o(1))} (1-e^{-c(1+o(1))})^{r} \nonumber \\
&\stackrel{\eqref{trick}}{\ge} \left(1-\frac{\delta^{i}}{N}\right)^{r(1+o(1))} (1-e^{-c})^{r(1+o(1))} \nonumber \\
&\ge \left(1-\frac{\delta^{i}}{N} - e^{-c}\right)^{r(1+o(1))} \nonumber \\
&\stackrel{\eqref{gammadef}}{=}\gamma_{smd}^{r(1+o(1))}.
 \end{align}
We used the fact that when $h_1(N) \rightarrow 0$, we have $\frac{1-e^{-c(1+h_1(N)})}{1-e^{-c}} \rightarrow 1$, which implies that there exists $h_2(N)\rightarrow 0 $ with $h_2(N)<\frac{1-e^{-c(1+h_1(N))}}{1-e^{-c}}-1$. Then, taking $h_3(N)=\frac{\log(1+h_2(N))}{\log(1-e^{-c})} \rightarrow 0$ we have
\begin{align}
\label{trick}
\frac{(1-e^{-c(1+h_1(N))})^{r}}{(1-e^{-c})^{r(1+h_3(N))}} &=\left( \frac{1-e^{-c(1+h_1(N)})}{1-e^{-c}}  (1-e^{-c})^{-h_3(N)}   \right)^r \nonumber\\
 &> ((1+h_2(N)) (1-e^{-c})^{-h_3(N)})^r=1.
\end{align}

(ii) In the case $ e^{-c} \ge 1-e^{-c}-\frac{\delta^i}{N}$, since by equation \eqref{jstar} we have $l^\star=i+2$,
\begin{align}
\label{caseip2}
\P(Z_v=1)&= \prod\limits_{l=1}^{l^\star-1}(1-p)^{r\delta^{l-1}\left(1+C\zeta\right)}\nonumber \\
 &\ge e^{-p(1+o(1))r\delta^{i}}\nonumber \\
 &\ge (e^{-c})^{r(1+o(1))}\nonumber \\
 &\stackrel{\eqref{gammadef}}{=}\gamma_{smd}^{r(1+o(1))}.
 \end{align}
 
Combining equations \eqref{caseip1} and \eqref{caseip2} we can deduce that
\begin{align}
\label{survivalp2}
\P(Z_v=1) \ge \gamma_{smd}^{r(1+o(1))}.
\end{align} 
 
Let $Z=\sum_{v \in V\setminus R} Z_v$ be the number of survivors in the SSP. By equation \eqref{survivalp2} we have 
\begin{equation}
\label{SMDEZ}
\E[Z]\ge(N-r)\gamma_{smd}^{r(1+o(1))}. 
\end{equation}

We finish with the computation similarly to equation \eqref{SQC_eq3} in Section \ref{SQC_lower}. The $o(1)$ term will be swallowed by the $\epsilon$ term in $r$. In particular, we will need the inequality
\begin{equation}
\label{epeta}
(1+o(1))(\eta-\epsilon)<(\eta-\epsilon/2),
\end{equation}
which holds because $\eta$ is upper bounded by one and we can choose an $\epsilon$ that tends to zero slower than the function hidden in the o(1) term. Then, putting it all together,

\begin{align}
\label{SQC_redo}
(\P(\mathcal{R})- \P(\overline{\mathcal{E}}))\P(\mathcal{E}) & \stackrel{\eqref{exres}}{\le} N^r\P(Y_2=0) \nonumber \\
&\stackrel{\eqref{EZ}}{\le} N^r 2e^{\frac43} e^{-\frac{\E[Z]}{3}}\nonumber \\
&\stackrel{\eqref{SMDEZ}}{=}2e^{\frac43}\exp\left(r\log(N)-\frac13(N-r)\gamma_{smd}^{r(1+o(1))}\right)\nonumber \\
&= 2e^{\frac43}\exp\left(r\log(N)-\frac13(N-r)\gamma_{smd}^{(1+o(1))(\eta-\epsilon) \frac{\log(N)}{\log(1/\gamma_{smd})}}\right)\nonumber \\
&\stackrel{\eqref{epeta}}{\le} 2e^{\frac43}\exp\left(r\log(N)-\frac13(N-r)\gamma_{smd}^{(\eta-\epsilon/2) \frac{\log(N)}{\log(1/\gamma_{smd})}}\right)\nonumber \\
&= 2e^{\frac43}\exp\left(r\left(\log(N)+\frac13 N^{-(\eta-\epsilon/2)}\right) -\frac13N^{1-(\eta-\epsilon/2) }\right)\nonumber \\
&\le 2e^{\frac43}\exp\left(\frac{(\eta-\epsilon)2\log^2(N)}{\log(1/\gamma_{smd})}-\frac13N^{1-(1+\log_N\log(1/\gamma_{smd}))+\epsilon/2}\right)\nonumber \\
&\le 2e^{\frac43}\exp\left(\frac{\log^2(N)-\frac13N^{\epsilon/2}}{\log(1/\gamma_{smd})}\right) \rightarrow 0
\end{align}
since $\frac{1}{\log(1/\gamma_{smd})}>1$ and  $\log^2(N)-\frac13N^{\epsilon}\rightarrow -\infty$ as long as $\epsilon \gg \frac{\log\log(N)}{\log(N)}$.  Finally, since $\P(\mathcal{E}) \rightarrow 1$, we have that there exists no resolving set a.a.s.
\end{proof}
%%%%%%%%%%%%%%%%%%%%%%%%%%%%%%%%%%%%%%%
\subsection{Proof of Theorem \ref{main_thm} for $p=o(1)$, Upper Bound}

If $c\rightarrow \infty$, the Theorem follows directly by $\mathrm{SMD} \le \mathrm{MD}$ and the results of \cite{bollobas2012metric}. For the remainder of this proof, we assume $c=\Theta(1)$.

It would be ideal to reduce this proof to the SQC proof, similarly to the lower bound. However, in case $p=\Theta(N^{-\frac{i}{i+1}})$ with $i>0$, the same approach as in Section \ref{SQC_upper} does not work. The events that two nodes $v$ and $w$ separate a set $W$ are not independent anymore and we cannot show that every set $W$ has an $f$-separator. Fortunately, we do not need such a powerful result for Theorem \ref{main_thm}; it is enough to prove the existence of $f$-separators for the subsets that can be candidates set in MAX-GAIN. In order to know which subsets can candidate sets, we expose most of the graph $G$, except that we reserve a small set $F$ of $c_\gamma\log^2(N)$ nodes with $c_\gamma=2/\log(1/\gamma_{smd})$, which we keep completely unexposed. The advantage of this is twofold. Now we can have a very good idea about which sets we need to separate, and we still have a large enough set of unexposed nodes that can independently separate the potential candidate sets (conditioned on the expansion properties of the exposed graph).

We have to make the claim that we have ``a very good idea about which sets we need to separate'' rigourous. Also, we must develop tools that allow us to reason about distances in the graph even when a small subset of the nodes are kept unexposed. We start showing that the unexposed nodes are a.a.s. far from each other.
 
%%%%%%%%%%%%%%%%%%%%%%%% NO HARM
\begin{lem}
\label{ip1o2}
In $G\sim\mathcal{G}(N,p)$, for a randomly selected $F\subset V$ with size $|F|=c_\gamma\log^2(N)$ with $c_\gamma= \frac{2}{\log(1/\gamma_{smd})}$ and for any two nodes $v,w\in F$, we have $d(v,w)\in \{i+1, i+2\}$ a.a.s., with $i$ given in Definition \ref{parameters}.
\end{lem}

\begin{proof}[Proof of Lemma \ref{ip1o2}]
We can sample the $c_\gamma\log^2(N)$ nodes one by one. Each time we sample one, let us also expose its $i$-neighborhood as done in the proof of Lemma \ref{fraction}. When we already sampled the first $j<\log^2(N)$ nodes, there are at most $c_\gamma\log^2(N)\delta^i(1+O(\zeta))$ nodes exposed (by Corollary \ref{cor_S}), so the probability that we select an unexposed node is at least $1-\frac{c_\gamma\log^2(N)\delta^i(1+O(\zeta))}{N}$. The probability that we always select an unexposed node is at least 
\begin{equation}
\left(1-\frac{c_\gamma\log^2(N)\delta^i(1+O(\zeta))}{N}\right)^{c_\gamma\log^2(N)}\rightarrow 1,
\end{equation}
because $\delta \ge \log^5(N)$ and $cN=\delta^{i+1}$ implies $\frac{\log^2(N)\delta^{i}}{N} = \frac{c\log^2(N)}{\delta} \le \frac{1}{\log^3(N)}$ for $N$ large enough. Since unexposed nodes are always at least $i+1$ distance away from all previous nodes and not more than $i+2$ by Corollary \ref{cor_diam}, this completes the proof.
\end{proof}

We now introduce the necessary definitions to reason about distances in $G$ with a small subset of the nodes $F$ unexposed.

\begin{defn}[exposed graph]
\label{exposed_graph}
Let $V'=V\setminus F$, let $G'$ be the subgraph of $G$ restricted to nodes $V'$. Let $N'=|V'|,\delta_2, c', i', \gamma_{smd}', \zeta'$ be the parameters defined in Definition \ref{parameters} and equation \eqref{gammadef} for graph $G'$. Let $d'(v,u)$ be the length of shortest path between nodes $v \in F$ and $u \in V'$. For $v \in V'$, $\mathcal{S}_{G'}(v,l)$ is defined in Definition~\ref{levelsetdef} for graph $G'$. For $v \in F$ let us extend the definition of $\mathcal{S}_{G'}(v,l)$ using the distance function $d'$ instead of $d$.
\end{defn}

In the rest of the section we will mainly use parameters $N',\delta_2, c', i', \gamma_{smd}',\zeta'$. We must keep in mind, that to prove Theorem \ref{main_thm} we must show 
$$\mathrm{SMDP} \le (1+o(1)) \log(N)/\log(1/\gamma_{smd}).$$ Fortunately, we can show that $\gamma_{smd}'=\gamma_{smd}^{1+o(1)}$, which means that proving $SMDP\le (1+o(1))\log(N')/\log(1/\gamma'_{smd})$ is enough for the theorem to hold. The following lemma will also show that $i'=i$ (consequently, we will not use the notation $i'$ after the following lemma).

\begin{lem}
\label{smdGLemma}
With the definitions given in Definition \ref{exposed_graph} we have
\begin{equation}
\label{smdG}
i'=i \quad \text{and} \quad \gamma_{smd}'=\gamma_{smd}^{1+o(1)}.
\end{equation}
\end{lem}

\begin{proof}[Proof of Lemma \ref{smdGLemma}]
First, we show that for any constant $k\ge1$, we have
\begin{equation}
\label{sandwitch}
\delta^{k}=(Np)^k\ge (N-c_\gamma\log^2(N))^kp^k = \delta_2^{k} \ge \delta^{k}\left(1-\frac{kc_\gamma\log^2(N)}{N}\right).
\end{equation}

Only the last inequality is not trivial. For the last inequality we note that since $x^k$ is a convex function for $k>1$,
\begin{equation}
\delta^{k}-\delta_2^{k}=(Np)^k- (N-c_\gamma\log^2(N))^kp^k \le \left(c_\gamma\log^2(N)kN^{k-1}\right)p^k = \delta^k\frac{kc_\gamma\log^2(N)}{N}.
\end{equation}

Equation \eqref{sandwitch} implies that $\delta^i/N=\Theta(1)$ if and only if  $\delta_2^i/N'=\Theta(1)$, hence $i'=i$. Moreover,
\begin{equation}
\label{sandwitchc}
\frac{N}{N'}c=\frac{\delta^{i+1}}{N'}>c'=\frac{\delta_2^{i+1}}{N'}\ge \frac{\delta^{i+1}}{N}\left(1-\frac{ic_\gamma\log^2(N)}{N}\right)\ge c\left(1-\frac{c_\gamma\log^3(N)}{N}\right),
\end{equation}
since $i\le \log(N)$. Equations \eqref{sandwitchc} and \eqref{gammadef} and imply equation \eqref{smdG}.
\end{proof}

Definition \ref{exposed_graph} also introduces the distance function $d'(v,u)$ on $G'$ (and one extra node from $F$). This function will be useful for us, first because it does not use any edge incident to $F\setminus \{v\}$, and therefore can be evaluated even if none of the edges incident to $F \setminus \{v\}$ are exposed, and second because we can prove that with high probability it is the same as the true distance. For the rest of this section, all expectations and probabilities are conditioned on the event that the expansion properties hold in the exposed graph $G'$.

\begin{lem}
\label{distclose}
In $\mathcal{G}(N,p)$, for a randomly selected $F\subset V$ with size $|F|=c_\gamma\log^2(N)$ with $c_\gamma= \frac{2}{\log(1/\gamma_{smd}')}$ for any two nodes $v \in F$ and $u \in V'$, we have $d'(v,u)=d(v,u)$ a.a.s.
\end{lem}

\begin{proof}[Proof of Lemma \ref{distclose}]
Since both $d(v,u)$ and $d'(v,u)$ represent distances, and the only difference is that the former is the distance in $G$ and latter is the distance in a subgraph of $G$, we must have $d(v,u) \le d'(v,u)$. Suppose that for some $v\in F ,u \in V'$ we have $d(v,u) < d'(v,u)$. This can happen only if there exists $w \in F \setminus \{v\}$ for which $d(v,w) + d(w,u) < d'(v,u)$. By Lemma \ref{ip1o2} we have that $d(v,w) \ge i+1$ a.a.s., and we also know that $d(u,w)\ge 1$ since $w\ne u$. If we could also show $d'(v,u) \le i+2$ a.a.s., the contradiction given by the inequality 
$$1+(i+1)\le d(v,w) + d(w,u) < d'(v,u) \le i+2$$
would prove that no such $w$ can exist a.a.s. Unfortunately, we cannot simply apply Corollary \ref{cor_diam} to show $d'(v,u) \le i+2$ a.a.s., since $d'(v,u)$ could be larger than $d(v,u)$. However, a simple analysis conditioning on the expansion properties suffices on the set $V'$ suffices. Indeed,
\small{
\begin{align}
\P(\exists v \in F, u\in V' \text{ with } d'(v,u)>i+2) & \le \sum\limits_{u \in V'} \P(\exists v \in F \text{ with } d'(v,u)>i+2)) \nonumber \\
& = \sum\limits_{u \in V'} \left(1- \P(\forall v\in F, v \in \bigcup_{l=1}^{i+2}\mathcal{S}_{G'}(u, l)) \right) \nonumber \\
& = \sum\limits_{u \in V'} \left(1-\left(1-(1-p)^{|\bigcup_{l=1}^{i+1}\mathcal{S}_{G'}(u, l)|}\right)^{|F|} \right) \nonumber \\
& \stackrel{\eqref{aasSip1}}{=} |V'| \left(1-\left(1-e^{-p(1-e^{-c'}+o(1))\delta_2^{i+1}}\right)^{c_\gamma\log^2(N)} \right).
\end{align}}
\normalsize{
By $1-x = e^{-x(1+O(1))}$ for $x=o(1)$, and $p\delta_2^{i+1}=c'\delta_2$ we proceed to}
\small{
\begin{align}
\P(\exists v \in F, u\in V' \text{ with } d'(v,u)>i+2) & \le N \left(1-e^{-e^{-c'\delta_2(1-e^{-c'}+o(1))}(1+o(1))c_\gamma\log^2(N)} \right) \nonumber \\
& \le N e^{-c'\delta_2(1-e^{-c'}+o(1))}(1+o(1))c_\gamma\log^2(N) \nonumber \\
& \le N^{1-c'\log^3(N)(1-e^{-c'}+o(1))} (1+o(1))c_\gamma\log^2(N) \nonumber \\
& \rightarrow 0,
\end{align}}
\normalsize{where in the last inequality we used $\delta_2=N'p\gg N' \log^5(N)/N \gg \log^4(N)$.}
\end{proof}

\begin{figure}[h]
\begin{center}
  \includegraphics[width=0.7\textwidth]{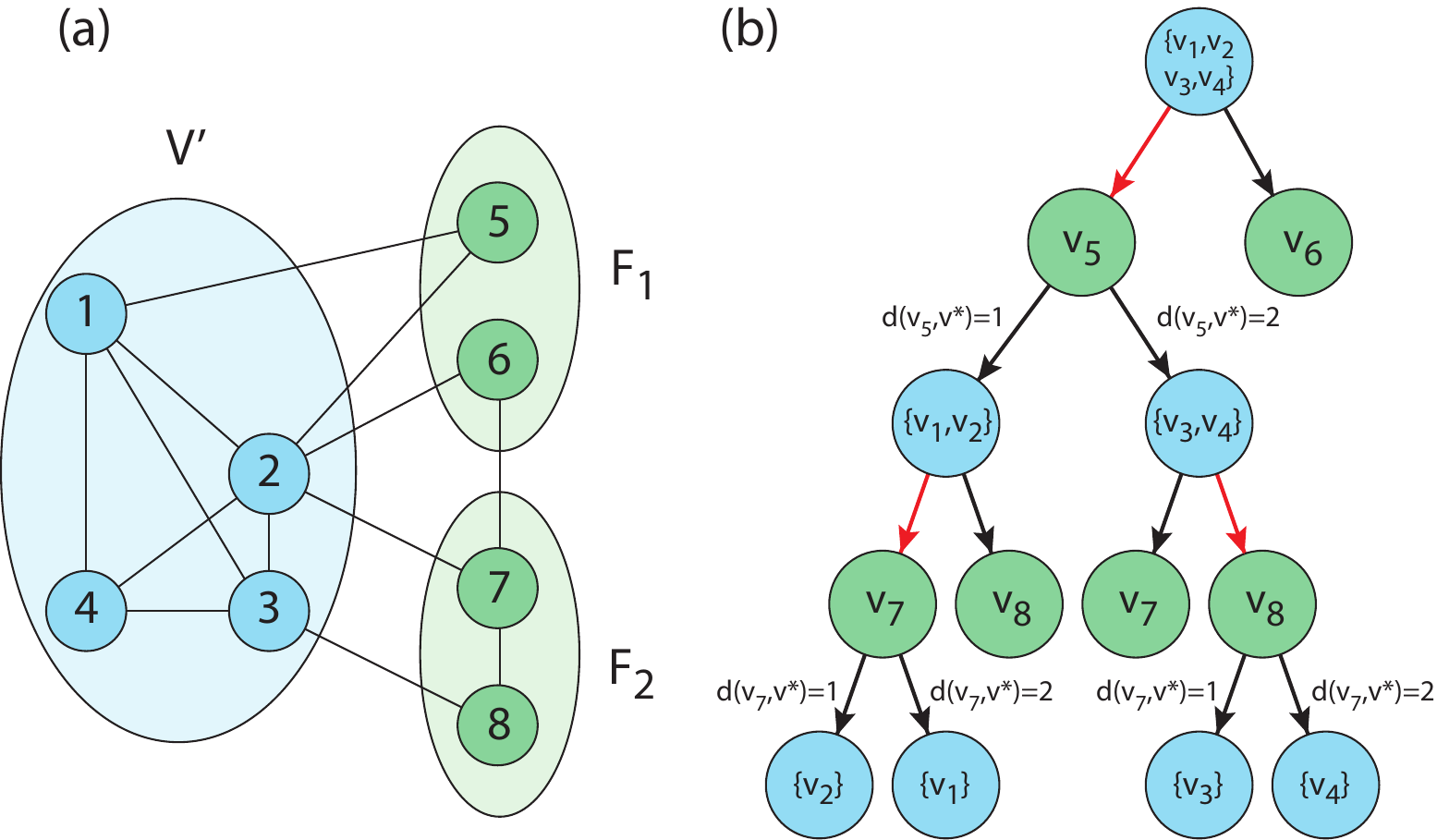}
  \caption{(a) A small example graph with $V'$, $F_1$ and $F_2$. (b) The ``game plan'' corresponding to the graph in (a). The $j^{th}$ blue layer ($j\ge0$) shows the potential pseudo-candidate sets we might encounter in step $j$. Arrows exiting blue nodes point to the potential of queries $F_j$ (green nodes on the $j^{th}$ level). The red arrow marks the query we picked. Arrows exiting green nodes correspond to the potential observations provided by Player 2 in the actual game. Each scenario ends when the potential candidate set has exactly one element. Tables \ref{table:T'} and \ref{table:R} show which sets $\mathcal{T}'_{j,\tilde{v}}$ and $R_{j,\tilde{v}}$ correspond to each step of this ``game plan''.}
  \label{game_tree}

\end{center}
\end{figure}

We have proved that we may use $d'(v,u)$ instead of $d(v,u)$, but we still do not know which sets need to be separated. To determine which sets we need to separate, we will simulate the game on the exposed graph. In the $j^{th}$ step, we assume that we have a candidate set, we expose a subset of the reserved nodes $F_j \subset F$ of size $\log(N)$, and we select the best reducer $v$ from only the nodes $F_j$ (unless the candidate set is small enough to query the whole set). Then, we consider all possible answers we could get if we selected $v$ as a query and we continue our simulation for each possible scenario (Figure \ref{game_tree} (b)). This analysis is different from the proof in Section \ref{SQC_upper} where we first proved a structural result for every subset and then proved that the MAX-GAIN algorithm finds the target. This time, the structural argument and the simulation of the algorithm will be intertwined. We simulate all possible scenarios of MAX-GAIN before actually taking observations from Player 2, and we construct function $g$ from Definition \ref{def_smd} to form a ``game plan'' that we can follow later in real-time. To implement this analysis, we need to slightly extend our definitions. 

From now on we will index our set of queries as $R_{j,\tilde{v}}$, since as we must prepare for observations for any target $\tilde{v}$. We now have the property $|R_{j,\tilde{v}}|=j$ and $R_{j,\tilde{v}} \subset R_{j+1,\tilde{v}}$, for all $\tilde{v}\in V$. We will also define a new version of candidate targets that are indexed by $\tilde{v}$, and which in addition uses the new distance function that we defined above. In this new notion of candidate targets we assume that the target is in the exposed graph $V'$ (the case when the target is in $F$ will be handled at the end of the proof).

\begin{defn}
Given a graph $G=(V,E)$ with unexposed nodes $F$ and queries $R_{j,\tilde{v}}$, the set of \textit{pseudo-candidate targets} 
$$\mathcal{T}'_{j,\tilde{v}}=\{ v \in V' \mid d'(w, v) = d'(w,\tilde{v}) \text{ for all } w \in R_{j,\tilde{v}} \}.$$ 
\end{defn}

\begin{rem}
\label{rem_partition}
Notice that $\tilde{v} \in \mathcal{T}'_{j,\tilde{v}}$ always holds. Also the sets $\mathcal{T}'_{j,\tilde{v}}$ define a partition on $V'$, that is for any $\tilde{v}\in V'$, 

\begin{itemize}
\item
$w \in \mathcal{T}'_{j,\tilde{v}} \Rightarrow \mathcal{T}'_{j,\tilde{v}}=\mathcal{T}'_{j,w}$ 
\item
 $w \not\in \mathcal{T}'_{j,\tilde{v}}\Rightarrow\mathcal{T}'_{j,\tilde{v}} \cap\mathcal{T}'_{j,w} =\emptyset$. 
 \end{itemize}
This can be seen by an inductive argument. For $j=0$, all sets $\mathcal{T}'_{0,\tilde{v}}$ coincide with $V'$, as $R_{0,\tilde{v}}=\emptyset$. For step $j+1$, each equivalence class at step $j$ is partitioned further by the new query.
\end{rem}

We must also define an analogous notion to extend Definition \ref{f-sep-graphs} of $f$-separators.

\begin{defn}[$f$-pseudo-separator]
Let $f(n)\in \mathbb{N}\rightarrow\mathbb{R}^+$ be a function. A set of nodes $W\subseteq V'$, $|W|=n$ has an $f$\textit{-pseudo-separator} if there is a node $w \in F$ such that 
\begin{equation}
\max\limits_{l \in \mathbb{N}}|W\cap \mathcal{S}_{G'}(w,l) |\le n\gamma_{smd}'+f(n).
\end{equation}
\end{defn}

\begin{algorithm}
\caption{Simulating all scenarios of MAX-GAIN}
\label{simulation}
\begin{enumerate}
\item We arbitrarily select $\log(N)$ disjoint sets $F_j \subset V$ of size $\log(N)$ and we let $F=\cup F_j$. We expose all edges of $V'=V \setminus F$. 
\item In step $j\ge0$, we expose the edges of nodes of $F_j$. For each $\mathcal{T}'_{j,\tilde{v}}$ we pick the best reducer $s_{j,\tilde{v}}\in F_j$ (possibly a different reducer one for each $\mathcal{T}'_{j,\tilde{v}}$) and add $s_{j,\tilde{v}}$ as a query to $R_{j,\tilde{v}}$. In the analysis, we prove that there always exists an $f$-separator (we define $f$ later), so the new query will also be a $(|\mathcal{T}'_{j,\tilde{v}}|\gamma_{smd}'+f(|\mathcal{T}'_{j,\tilde{v}}|)$-reducer for $\mathcal{T}'_{j,\tilde{v}}$. Selecting it produces the new sets $\mathcal{T}'_{j+1,\tilde{v}}$.

\item When a set $\mathcal{T}'_{j,\tilde{v}}$ reaches size $o(\log(N)/\log(1/\gamma_{smd}))$, we query the entire set (see the proof of Lemma \ref{prog}, for the base case).
\end{enumerate}
\end{algorithm}

\begin{table}[h]
\begin{center}
\begin{tabular}{ |c|c|c|c|c| }
  \hline
   $j$ & $\mathcal{T}'_{j,v_1}$ & $\mathcal{T}'_{j,v_2}$ &$\mathcal{T}'_{j,v_3}$ & $\mathcal{T}'_{j,v_4}$ \\  \hline
  0 & $V'$ & $V'$ & $V'$ & $V'$  \\  \hline 
  1 & $\{v_1,v_2\}$ & $\{v_1,v_2\}$  & $\{v_3,v_4\}$ & $\{v_3,v_4\}$ \\  \hline 
  2 & $\{v_1\}$  & $\{v_2\}$  &$\{v_3\}$ & $\{v_4\}$ \\  \hline 
\end{tabular}
\end{center}
  \caption{The pseudo-candidate sets corresponding to each $\tilde{v}$ and $j$ from the example in Figure \ref{game_tree}.}
%\end{table}
\label{table:T'}

%\begin{table}[h!]
\begin{center}
\begin{tabular}{ |c|c|c|c|c| }
  \hline
  $j$  &  $R_{j,v_1}$& $R_{j,v_2}$& $R_{j,v_3}$& $R_{j,v_4}$   \\  \hline
  0 &$\emptyset$&$\emptyset$&$\emptyset$&$\emptyset$ \\  \hline 
  1 &$\{v_5\}$&$\{v_5\}$  &$\{v_5\}$& $\{v_5\}$ \\  \hline 
  2 &$\{v_5,v_7\}$  & $\{v_5,v_7\}$  &$\{v_5,v_8\}$ & $\{v_5,v_8\}$ \\  \hline 
\end{tabular}
\end{center}
  \caption{The query sets corresponding to each $\tilde{v}$ and $j$ from the example in Figure \ref{game_tree}.}
  \label{table:R}
\end{table}

With these definitions, the simulation of MAX-GAIN is defined in Algorithm \ref{simulation}. We must show that while constructing our ``game plan'', it is in fact possible with probability tending to one to select an $f$-pseudo-separator $s_{j,\tilde{v}}$ in each step of the algorithm only from the $F_j$ (we define $f$ later). Then Lemma \ref{distclose} implies that these $f$-pseudo-separators for the pseudo-candidate sets are in fact $f$-separators for the true candidate sets. This will allow us to use Lemma \ref{prog}.

 %%%%%
\begin{lem}
\label{eval_reserved}
Let $n=|\mathcal{T}'_{j,\tilde{v}}|$ be the cardinality of the pseudo-candidate target set in step $j$ of the game plan for target~$\tilde{v}$. Let $v \in F_j$, let $X_{vw}$ be the indicator of the event $v\in \mathcal{S}_{G'}(w,i+2)$ for $w \in \mathcal{T}'_{j,\tilde{v}}$, and let $X_v=\sum_{w \in \mathcal{T}'_{j,\tilde{v}}} X_{vw}=|\mathcal{S}_{G'}(v,i+2) \cap \mathcal{T}'_{j,\tilde{v}}|$. Then,
\begin{equation}
\E[X_v]=n(e^{-c'} + O(\zeta')).
\end{equation}
\end{lem}

\begin{proof}[Proof of Lemma \ref{eval_reserved}]
This is an analogous result to equation \eqref{condExp} in Lemma \ref{fraction}, except that now most of the graph is exposed. Consider $w\in \mathcal{T}'_{j,\tilde{v}}$, then

\begin{align}
\P(X_{vw}=1) &= (1-p)^{|\cup_{j=1}^i\mathcal{S}_{G'}(w,j)|}\nonumber\\
&= e^{(-p+O(p^2))(\delta_2^{i}(1+O(\zeta')))}\nonumber\\
&= e^{-c'}(1+O(\zeta')))\nonumber\\
&= e^{-c'} + O(\zeta').
\end{align}
The result on the expectation follows immediately .
\end{proof}

The previous lemma only covered the expectation. Now we will establish a result on concentration.

\begin{lem}
\label{concentration_reserved}
Let $v, X_v, X_{vw}, n$ be defined as in Lemma \ref{eval_reserved}, and let $\omega$ be a function tending slowly to infinity. Then,
\begin{equation}
\P\left(|X_v- ne^{-c'}| > \frac{n}{\omega\log(n)} \right)\rightarrow 0
\end{equation}
as $N\rightarrow\infty$ independently of $n$.
\end{lem}

\begin{proof}[Proof of Lemma \ref{concentration_reserved}]
By Lemma \ref{eval_reserved} and Chebyshev's inequality,
\begin{align}
\label{concentration_reserved_eq}
&\P\left(|X_v- ne^{-c'} | > \frac{n}{\omega\log(n)} \right)\nonumber \\
&\le \P\left(|X_v- n(e^{-c'} + O(\zeta')) | > \frac{n}{2\omega\log(n)} \right) \nonumber \\ 
&=\P\left(|X_v- \E[X_v] | > \frac{n}{2\omega\log(n)}\right) < \frac{4\omega^2\log^2(n)\Var[X_v]}{ n^2}.
\end{align}
To compute $\Var[X_v]$ we will need
\begin{align}
\E[X_{vw}X_{vx}]&= \P(d(v,w)= i+2 \text{ and } d(v,x)= i+2) \nonumber \\
&=(1-p)^{|\cup_{j=1}^i\mathcal{S}_{G'}(w,x,j)|}\nonumber \\
&\stackrel{\eqref{expansioneq2}}{=} e^{(-p+O(p^2))(2\delta_2^i+O(\zeta'))}\nonumber\\
&= e^{-2c'}(1+O(\zeta'))\nonumber\\
&= e^{-2c'} + O(\zeta').
\end{align}
Then,
\begin{align}
\Var[X_v]&=\E[X_v^2]-\E^2[X_v]\nonumber \\
&= \sum_{w \in \mathcal{T}'_{j,\tilde{v}}} \left( \E[X_{vw}^2] - \E^2[X_{vw}] \right) \nonumber\\
& \qquad\qquad+ \sum_{w,x}  \left(\E[X_{vw}X_{vx}] - \E[X_{vw}]\E[X_{vx}]\right)\nonumber\\
&\le \E[X_v]+ n^2(e^{-2c'}-e^{-c'}e^{-c'}+ O(\zeta'))\nonumber\\
&=ne^{-c'}+n^2O(\zeta'). \label{eq:var_calc}
\end{align}

Recall that we assumed $c=\Theta(1)$ at beginning of the section, and that in equation~\eqref{sandwitchc} we proved that $c'=\Theta(1)$ must hold as well. Consequently, we have
$$\frac{\delta_2^{i}}{N'}=\Theta\left(\frac{1}{\delta_2}\right) = o\left(\sqrt{\frac{\log(N')}{\delta_2}} \right),$$
which implies that we can choose $\zeta'=\sqrt{\log(N')/\delta_2}$ in the version of equation \eqref{zeta1} where we replace the parameters $ \zeta,\delta, N$ by the parameters $ \zeta’,\delta_2, N’$. Then, by the assumption $pN\gg \log^5(N)$ in the statement of Theorem~\ref{main_thm}, we have
\begin{equation}
\label{eq:zeta'bound}
\zeta' = \sqrt{\frac{\log(N')}{\delta_2}} = o\left(\sqrt{\frac{N\log(N')}{N' \log^5(N)}} \right)=o\left(\frac{1}{\log^2(N)} \right).
\end{equation}
Finally, substituting equation \eqref{eq:var_calc} into \eqref{concentration_reserved_eq} yields
\begin{align}
\P\left(|X_v- ne^{-c'} | > \frac{n}{\omega\log(n)} \right) &< \frac{4\omega^2\log^2(n)\Var[X_v]}{ n^2} \nonumber \\
&\stackrel{\eqref{eq:var_calc}}{=}\frac{4\omega^2\log^2(n)ne^{-c'}}{ n^2} + \frac{4\omega^2\log^2(n)n^2O(\zeta')}{n^2} \nonumber \\
&\stackrel{\eqref{eq:zeta'bound}}{=} o(1)+ o\left( \frac{\omega^2\log^2(n)}{\log^2(N)} \right),
\end{align}
which proves the desired result since $N\ge n$ and $\omega$ tends to infinity very slowly.
\end{proof}

\begin{lem}
\label{not_stuck}
Let $Z$ be the indicator variable that we cannot select an $f$-pseudo-separator in some step of the simulation with $f(n)=2n/(\omega\log(n))$. Then $\P(Z) \rightarrow 0$.
\end{lem}

\begin{proof}[Proof of Lemma \ref{not_stuck}]
Let $Z_j$ be the indicator variable that we cannot select an $f$-pseudo-separator in the $j^{th}$ step of the simulation. Let us fix $j$. Let $Y_{j,\tilde{v}}$ be the indicator variable that we cannot find an $f$-pseudo-separator for the pseudo-candidate set $\mathcal{T}'_{j,\tilde{v}}$. Since by Remark \ref{rem_partition} the pseudo-candidate sets partition $V'$, some (for $j=1$ all) of the $Y_{j,\tilde{v}}$ can be identical, but this will not matter as in the end we will apply a union bound. Similarly to the proof of the SQC upper bound, finding an $f$-pseudo-separator is equivalent to finding an $X_v$ close to its expectation. Indeed, $|X_v-ne^{-c'}| \le f(n)/2$ implies
\begin{align}
 X_v& \le ne^{-c'}+ \frac{f(n)}{2} \le n\gamma_{smd}' +\frac{f(n)}{2} <n\gamma_{smd}' +f(n) 
\end{align}
and
\begin{align}
 n-X_v & \le n(1-e^{-c'})+ \frac{f(n)}{2} \le n\gamma_{smd}' + n\frac{\delta_2^i}{N'} +\frac{f(n)}{2} < n\gamma_{smd}' + f(n),
\end{align}
because, as we saw in equation \eqref{eq:zeta'bound}, we can choose $f(n)$ so that $\delta_2^i/N' <\zeta'<1/\log^2(N)<f(n)/(2n)$.  Thus, $v$ is an $f$-pseudo-separator. The non-existence of an $f$-pseudo-separator implies the non-existence of an $X_v$ close to its expectation, which means
\begin{equation}
\P(Y_{j,\tilde{v}})\le \P\left(|X_v- ne^{-c'}| > \frac{f(n)}{2} \quad\forall v\in F_j \right).
\end{equation}
Let us choose $N$ large enough such that for $v \in F_j$
\begin{equation}
\P\left(|X_v- ne^{-c'}| > \frac{f(n)}{2} \right)<e^{-2}
\end{equation}
(which can be done for any constant by Lemma \ref{concentration_reserved} since $f(n)/2 = n/(\omega\log(n))$). Then,
\begin{equation}
\P(Y_{j,\tilde{v}})\le e^{-2|F_j|}=N^{-2}.
\end{equation}
By union bound, since in every step we have at most $N'<N$ pseudo-candidate sets $\mathcal{T}'_{j,\tilde{v}}$ to separate,
\begin{equation}
\P(Z_j)=\P\left(\bigcup_{\tilde{v} \in |V'|}Y_{j,\tilde{v}} \right)\le N\P(Y_{j,\tilde{v}})\leq N^{-1}.
\end{equation}

Finally, since we have $\frac{2\log{N}}{\log(1/\gamma_{smd})}$ sets $F_j$, another union bound shows that
\begin{equation}
\P(Z)=\P\left(\bigcup_{j=1}^{\frac{2\log{N}}{\log(1/\gamma_{smd})}}Z_j\right) \le \frac{2\log{N}}{N\log(1/\gamma_{smd})}=o(1).
\end{equation}
\end{proof}

Now we just need to put the pieces together to prove Theorem \ref{main_thm}.

\begin{proof}[Proof of Theorem \ref{main_thm}, upper bound]
We perform Algorithm \ref{simulation}, with the modification that besides $F$, we also reserve another set $F'$ with $c_\gamma\log^2\log^2(N)$ nodes. The modified algorithm runs in three steps.

(i) We run Algorithm \ref{simulation} on $V'=V \setminus (F \cup F')$. The additional set $F'$ slightly increases the size of the reserved nodes, but this $\log^2\log^2(N)$ term does not affect the analysis. Lemma \ref{not_stuck} ensures that the algorithm can find an $f$-pseudo-separator for all $F_j$ with probability tending to $1$. Lemma \ref{distclose} shows that the only candidate sets we might encounter in the MAX-GAIN algorithm are pseudo-candidate target sets in our game plan, and the $f$-pseudo-separator we found for the pseudo-candidate target sets are $f$-separators for the corresponding candidate target sets, unless the source was in the reserved nodes.

Since the $f$ we used in Lemma \ref{not_stuck} was $o(n/\log(n))$ we can apply Lemma \ref{prog}, which shows that in each possible scenario we simulate, we find the source in $$(1+o(1))\frac{\log(N')}{\log(1/\gamma_{smd}')}$$ steps. Therefore the number of steps we require is always is less than $$\frac{2\log(N)}{\log(1/\gamma_{smd})},$$ the number of sets $F_j$ we can use to find $f$-separators in Lemma \ref{not_stuck}. Thus, if the target was in $V'$, the algorithm will find it. By Lemma \ref{smdGLemma}, $\gamma_{smd}'=\gamma_{smd}^{1+o(1)}$, hence the number steps taken is upper bounded by the desired $$(1+o(1))\frac{\log(N)}{\log(1/\gamma_{smd})}$$ steps.

(ii) We repeat the argument with candidate set $F$ and reserved nodes $F'$.

(iii) Finally we query the entire $F'$. 

In this last two steps, we selected only $o(\log(N))$ extra queries, which does not change the leading term of our upper bound. We ensured that no matter whether the source is in $V\setminus (F \cup F')$, $F$ or $F'$, we will be able to find it in the desired number of steps with probability tending to 1. Recall that in all of our calculations in Lemmas \ref{distclose}-\ref{not_stuck} we conditioned on the event that the expansion properties hold in the exposed graph. Since the expansion properties also hold with high probability, the upper bound in Theorem \ref{main_thm} holds also without conditioning. 
\end{proof}

\section{Discussion}
\label{discussion}

In this paper, we proved tight asymptotic results for the SMD in $\mathcal{G}(N,p)$. We found that a.a.s., the ratio between the SMD and the MD is a constant as $N$ tends to infinity, and we conjecture that this constant is $1$ except for $(pN)^i=\Theta(N)$ for $i \in \mathbb{N}$, where the constant term is found explicitly and is smaller than $1$. On the one hand, considering the equivalence of binary search with adaptive and non-adaptive queries, it is interesting that there is any difference at all between the SMD and the MD. On the other hand, experimental results suggest that on other graph models (and especially real-world networks), the SMD is orders of magnitude smaller than the MD \cite{spinelli2017general}. Hence, the Erd\H{o}s-R\'enyi graphs are an intermediate regime, where the restriction on the queries does favor adaptive algorithms, but not by too much.

There are several open questions remaining. The lower and the upper bounds in Theorem \ref{main_thm} are a factor of $\eta$ apart; the same factor that appeared in the earlier work of \cite{bollobas2012metric}. We believe that a further study of the new notions introduced in this paper, the QC (which is essentially equivalent to the minimum cardinality of an identifying code) and the SQC, may help removing this gap. 

It would be interesting to study random graph models other than the $\mathcal{G}(N,p)$ model, where we expect the difference between the MD and the SMD to be significantly larger. Adding noise to the measurements would be another step towards more realistic scenarios, and in this case too, we expect the a larger difference between the MD and the SMD. The noise can come from faulty observers similarly to \cite{emamjomeh2016deterministic}, or the noise can be proportional to the distances observed which would model stochastic disease propagation in source localization \cite{lecomte2020noisy,Spinelli2016ObserverPF}.      

\acks

We thank Brunella Spinelli, Dieter Mitche, Pawe\l{} Pra\l{}at and J\'ulia Komj\'athy for their useful comments and suggestions.

The work presented in this paper was supported in part by the Swiss National Science Foundation under grant number 200021-182407.

%\appendix
%\input{implementation_challenges}

\bibliographystyle{APT}
\bibliography{literature}

\appendix
\section{Additional Proofs}
\label{proofs}

\subsection{Proof of Lemma \ref{prog}}
\label{proof_l2}

Let $T(n)$ denote the number of steps in which MAX-GAIN reduces the number of candidates from $n$ to 1, and let $C_N$ be the value (not depending on $n$) such that for all $n\ge C_N >0$ the condition
\begin{equation}
T(n) \le T(nq+f(n)) +1 \text{ with } f(n)=o\left(\frac{n}{\log(n)}\right)
\end{equation}
holds. Then we prove 
\begin{equation}
T(n)< \log_{\frac{1}{q}}(n) + \log\log(n)+C_{q,f} + C_N,
\end{equation}
where $C_{q,f}$ is a positive constant (it depends only on $q$ and $f$ but not $n$) computed implicitly at the end of the proof.

Proof by induction. Base case: if $n<C_{q,f} + C_N$ then $T(n)<C_{q,f} + C_N$ clearly holds as we can query each candidate. Induction step: Let now $n\ge C_{q,f} + C_N$ and we assume that for $M< n$ the induction hypothesis holds, that is

\begin{equation}
\label{ind.ass.}
T(M)< \log_{\frac{1}{q}}(M) + \log\log(M)+C_{q,f} + C_N.
\end{equation}
Then,

\begin{align}
\label{ind_step}
T&(n) \le T(nq+f(n)) + 1\nonumber \\
&\stackrel{\eqref{ind.ass.}} \le \log_{\frac{1}{q}}(nq+f(n))+\log(\log(nq+f(n))) +C_{q,f} + C_N + 1
\end{align}

For the induction hypothesis to hold we would like the last expression to be upper bounded by
$$\log_{\frac{1}{q}}(n) + \log(\log(n)) +C_{q,f} + C_N.$$
To compare these two quantities, we would like to transform $\log_{\frac{1}{q}}(nq+f(n))$. Using the fact that log is a concave function and by linearly approximating it at $n$,

\begin{align}
\log_{\frac{1}{q}}(nq+f(n)) &= \log_{\frac{1}{q}}(n+\frac{f(n)}{q})-1 \nonumber\\
&\le \log_{\frac{1}{q}}(n) + \frac{f(n)}{q\log(\frac{1}{q})n}-1
\end{align}

Plugging this into \eqref{ind_step} we get

\begin{equation}
T(n) \le  \log_{\frac{1}{q}}(n)+\frac{f(n)}{q\log(\frac{1}{q})n}+\log(\log(nq+f(n))) +C_{q,f} + C_N
\end{equation}
For the induction hypothesis to hold we need to show

\begin{align}
\log_{\frac{1}{q}}(n)+\frac{f(n)}{q\log(\frac{1}{q})n}+&\log(\log(nq+f(n))) +C_{q,f} + C_N \nonumber \\
& \le \log_{\frac{1}{q}}(n) + \log(\log(n)) +C_{q,f} + C_N,
\end{align}
which is equivalent to

\begin{equation}
\frac{f(n)}{q\log(\frac{1}{q})n}+\log(\log(nq+f(n))) \le \log(\log(n))
\end{equation}
Again, by the concavity of $\log(\log(n))$ we can use a linear approximation

\begin{equation} 
\log(\log(nq+f(n))) \le \log(\log(n)) + \frac{n-(nq+f(n))}{n\log(n)}
\end{equation}
So it is enough to show 

\begin{align}
\frac{f(n)}{q\log(\frac{1}{q})n} &\le \frac{n-(nq+f(n))}{n\log(n)}\nonumber\\
\frac{f(n)}{n} \left(\frac{1}{\log(\frac{1}{q})q}+\frac{1}{log(n)}\right) &\le \frac{1-q}{\log(n)}\nonumber\\
\frac{f(n)\log(n)}{n} &\le \left(\frac{1-q}{\log(\frac{1}{q})q}+\frac{1-q}{log(n)}\right)^{-1} 
\end{align}

Since the right hand side is bounded from below (for $n>0$) and $f(n)=o\left(\frac{n}{\log(n)}\right)$, this last inequality must hold for $n\ge C_{q,f}$, for some constant $C_{q,f}$ (depending only on $q$ and $f$ but not $n$).

To conclude the proof, we showed that for all $n \in \mathbb{N}$

\begin{equation}
T(n)< \log_{\frac{1}{q}}(n) + \log\log(n)+C_{q,f} + C_N.
\end{equation}

This in particular implies
\begin{equation}
T(N)< \log_{\frac{1}{q}}(N) + \log\log(N)+C_{q,f} + C_N = (1+o(1))\log_{\frac{1}{q}}(N)
\end{equation}

for $C_N=o\left(\log_{\frac{1}{q}}(N) \right)$.

\subsection{Proof of Lemma \ref{expansion}}

The proof follows the proof of Lemma 2 (i) in \cite{bollobas2012metric} until the very last step, the evaluation of the multiplicative error term. There, the authors use $i=O(\log(n)/\log\log(n))$ and $\sqrt{\omega}\le \log^2(N)\log\log(N)$ to get the asymptotic upper bound

\begin{align*}
\left(1+O\left(\frac{\delta}{N}\right)+O\left(\frac{1}{\sqrt{\omega}}\right) \right)&\prod\limits_{j=2}^i \left(1+O\left(\frac{\delta^j}{N}\right)+O\left(\frac{1}{\sqrt{\omega d^{j-1}}}\right) \right) \\
=& \left(1+O\left(\frac{\delta^i}{N}\right)+O\left(\frac{1}{\sqrt{\omega}}\right) \right)\prod\limits_{j=7}^{i-3}(1+O(\log^{-3}(N))) \\
=&\left(1+O\left(\frac{\delta^i}{N}\right)+O\left(\frac{1}{\sqrt{\omega}}\right) \right)(1+O(\log^{-2}(N))) \\
=& \left(1+O\left(\frac{\delta^i}{N}\right)+O\left(\frac{1}{\sqrt{\omega}}\right) \right) 
\end{align*}

However, the second upper bound on $\sqrt{\omega}$ is not necessary. Instead, we can write

\begin{align*}
\left(1+O\left(\frac{\delta}{N}\right)+O\left(\frac{1}{\sqrt{\omega}}\right) \right)&\prod\limits_{j=2}^i \left(1+O\left(\frac{\delta^j}{N}\right)+O\left(\frac{1}{\sqrt{\omega d^{j-1}}}\right) \right) \\
=& \left(1+O\left(\frac{\delta^i}{N}\right)+O\left(\frac{1}{\sqrt{\omega}}\right) \right)\prod\limits_{j=5}^{i-2}\left(1+O\left(\frac{1}{\delta^2}\right)\right) \\
=& \left(1+O\left(\frac{\delta^i}{N}\right)+O\left(\frac{1}{\sqrt{\omega}}\right) \right)\left(1+O\left(\frac{1}{\delta}\right)\right)\\
=& \left(1+O\left(\frac{\delta^i}{N}\right)+O\left(\frac{1}{\sqrt{\omega}}\right) \right) 
\end{align*}

where the second to last inequality holds because $i<\delta$ and $(1+O(1/\delta^2))^\delta=(1+O(1/\delta))$, and the last inequality holds because $1/\delta=O(\delta^i/N)$.

The condition $\sqrt{\omega}\le \log^2(N)\log\log(N)$ is not used anywhere else in the proof of Lemma 2 (i) in \cite{bollobas2012metric}, so we may remove this condition, which gives Lemma \ref{expansion} of this paper.

\end{document}